\documentclass[final,3p,times]{elsarticle}

\newenvironment{enumerate-i}{%

\begin{enumerate}%
}{%
\end{enumerate}}
\newenvironment{enumerate-a}{%

\begin{enumerate}%
}{%
\end{enumerate}}

\newcommand{\RE}{{\rm Re}}

\newcommand{\PE}{{\rm Pe}}

\newcommand{\RNum}[1]{\uppercase\expandafter{\romannumeral #1\relax}}
\newcommand{\bj}{\bm{j}}

\newcommand{\DA}{{\rm Da}}
\newcommand{\DAII}{{\rm Da}_{\rm \RNum{2}}}

\def\bx{{\bm{x}}}

\def\bv		{{ {\bm{v}} }}
\def\bV		{{ {\bm{V}} }}

\def\lob#1 {{ \left( #1 \right) }}
\def\lsb#1 {{ \left[ #1 \right] }}
\def\lcb#1 {{ \left\lbrace #1 \right\rbrace }}

\def\lab#1 {{ \left\langle #1 \right\rangle }}

\def\lcb#1 {{ \left\lbrace #1 \right\rbrace }}
\def\ddel#1 {{ \delta \left( #1 \right) }}
\def\bvs#1#2 {{ {\bv}_{#1}^{#2} }}

\newcommand{\bde}{\begin{description}}
\newcommand{\ede}{\end{description}}

\usepackage[bookmarks=false]{hyperref}
\hypersetup{
    colorlinks,
    citecolor=black,
    filecolor=black,
    linkcolor=black,
    urlcolor=black
}
\usepackage{bm}
\graphicspath{ {figures/} }
\usepackage{mathtools}
\newcommand{\grad}{\boldsymbol{\nabla} }
\newcommand{\bo}[1]{\bm{#1}}

\newcommand{\nx}{\grad_{\bm{x}}}
\newcommand{\ny}{\grad_{\bm{y}}}

\newcommand{\by}{\bm{y}}

\newcommand{\ofx}{\left( \bm{x} \right)}
\newcommand{\ofxt}{\left( \bm{x}, t \right)}
\newcommand{\ofxyt}{\left( \bm{x},\bm{y}, t \right)}
\newcommand{\ofy}{\left( \bm{y} \right)}

\renewcommand{\div}{\boldsymbol{\nabla} \cdot}

\newcommand{\ddt}[1]{\frac{\partial #1}{\partial t} }

\newcommand{\diff}{\mathcal{D}}
\newcommand{\eps}{\varepsilon}

\usepackage{lineno}
\usepackage[inline]{showlabels}

\usepackage{caption}
\usepackage{subcaption}

\journal{Advances in Water Resources}

\begin{document}

\begin{frontmatter}



\title{Macroscopic models for filtration and heterogeneous reactions in porous media}


\author[uon]{Federico Municchi\corref{cor1}}
\ead{federico.municchi@nottingham.ac.uk}
\cortext[cor1]{Corresponding author}
\author[uon]{Matteo Icardi} 
\address[uon]{School of Mathematical Sciences, University of Nottingham, Nottingham, UK}


\begin{abstract}
Derivation of macroscopic models for advection-diffusion processes in the presence of dominant heterogeneous (e.g., surface) reactions using homogenisation theory or volume averaging is often deemed unfeasible \citep{VALDESPARADA20112177,BATTIATO201118} due to the strong coupling between scales that
characterise such systems. In this work, we show how the upscaling can be carried out by applying and extending the methods presented in \citet{ALLAIRE2007523,Mauri_1991}. The approach relies on the decomposition of the microscale concentration into a reactive component, given by the eigenfunction of the  advection-diffusion operator, the associated eigenvalue which represents the macroscopic effective reaction rate, and a non-reactive component. The latter can be then upscaled with a two-scale asymptotic expansion and the final macroscopic equation is obtained for the leading order.
The same method can also be used to overcome another classical assumption, namely of non solenoidal velocity fields, such as the case of deposition of charged colloidal particles driven by electrostatic potential forces.
The whole upscaling procedure, which consists in solving three cell problems, is implemented for arbitrarily complex two- and three-dimensional periodic structures using the open-source finite volume library OpenFOAM\textsuperscript{\textregistered}.
We provide details on the implementation and test the methodology for two-dimensional periodic arrays of spheres, and we compare the results against fully resolved numerical simulations, demonstrating the accuracy and generality of the upscaling approach.
The effective velocity, dispersion and reaction coefficients are obtained for a wide range of P\'eclet and surface Damk{\"o}hler numbers, and for Coulomb-like forces to the grains.
Noticeably, all the effective transport parameters are significantly different from the non-reactive (conserved scalar) case, as the heterogeneity introduced by the reaction strongly affects the micro-scale profiles. 
\end{abstract}


\end{frontmatter}


\tableofcontents

\section{Introduction}
\label{S:intro}

Understanding transport in porous media is important for a wide range of chemical (for example drying of paper pulp or flow through catalysts), biological (membranes), and geological (groundwater remediation, $CO_2$ storage, nuclear waste management, etc.) applications \cite{Cheng2016}.   
However, transport phenomena are complicated by the broad range of excited space and time scales, which results in large (almost always unaffordable) computational requirements when simulating such systems at industrial or natural scales. Therefore, simulations describing all scales down to local heterogeneities (e.g., pore-scale) would be of little practical applicability, without a robust methodology that leads to reliable upscaled reduced models.

In the following, we will use the term \textit{'upscaling'} to describe the procedure leading to 'reduced/effective models' where the 'fast' (high wave number) components of the unknown fields are averaged out. We  also notice that often the terms 'spatial filtering' \cite{Radl2018} or 'coarse graining' \cite{kardar_2007} are also employed in literature with a similar meaning.  

A variety of methods can be employed to perform this upscaling procedure: asymptotic homogenisation \cite{Pavliotis2008,Davit2013,auriault1991} is a powerful and versatile tool for the upscaling of transport and reaction equations in porous and heterogeneous media.
For example, Taylor dispersion in porous media can be approached in this way \citep{auriault1995taylor}. Homogenisation has also been applied to rigorously derive models for heterogeneous (surface) reactions~\citep{BT2011a}, and colloidal particles deposition\cite{KREHEL2015209,RAY2012374}. These are, however, limited to the regime of slow reactions (or mass transfer) at the pore surfaces, leading therefore to a trivial limit where the reaction (scaled by the surface to volume ratio) simply appears in the macroscopic equation as a source term without affecting the other transport parameters. Another common limitation is the assumption of divergence-free velocities.

Although these limitations are usually well known and have been already discussed and partially addressed in the last decades \cite{ALLAIRE2007523,Mauri_1991,boccardo2018robust}, this is still mostly neglected when dealing with more realistic applications.
 In this work we show how the issues highlighted above can be overcome. We derive a homogenised equation using the technique developed by Allaire and Raphael \citep{ALLAIRE2007523} in the case of homogeneous surface reactions and we generalised it to include non-homogeneous boundary conditions and electrostatic forces.
    
    The objectives of this paper are threefold: i) re-deriving and discussing in a clear and extensive manner the seminal approach of \cite{ALLAIRE2007523,Mauri_1991} to make it more widely understandable by the porous media community; ii) extending the approach to cover relevant cases in applications, namely more general non-homogeneous boundary conditions, and electrostatic (or other external) forces; iii) implementing an algorithm to validate the approach and test it on two- and three-dimensional geometries. To ensure the full reproducibility of results we also provide the community with an open-source upscaling toolbox based on the OpenFOAM\textsuperscript{\textregistered} library.

    The following sections are structured as follows: after introducing the governing equations and physical problem in Section~\ref{S:model}, we present and discuss the approach of Allaire and Raphael in Section~\ref{S:homogeneous}, and develop novel extensions in Section~\ref{S:extensions}. In section~\ref{S:num} we detail the numerical algorithm, while Section~\ref{S:res} compares the method against fully resolved numerical simulation and present some results for the upscaled parameters obtained for a periodic arrays of spheres under different flow regimes. Eventually, we present conclusions and outlook on future directions in section~\ref{S:conc}.

\section{Mathematical model}
\label{S:model}

    Consider a porous medium occupying a region of space $\hat{\Omega}$ associated with a characteristic length $L$. We assume  $\hat{\Omega}$ to be composed of spatially repeated (periodic) unit cells $\hat{\mathcal{Y}}$ with characteristic length (periodicity) $\ell=\eps L \ll 1$. Each unit cell $\hat{\mathcal{Y}}$ is then given by the union $\hat{\mathcal{Y}} = \hat{\mathcal{Y}}_f \cup \hat{\mathcal{Y}}_s$, where $\hat{\mathcal{Y}}_f$ and $\hat{\mathcal{Y}}_s$ are the fluid and solid regions of $\hat{\mathcal{Y}}$ respectively, separated by an interface $\hat{\Gamma}$.
    Clearly, $\hat{\mathcal{Y}}_f$ is generally not simply connected, while $\hat{\mathcal{Y}}_s$ is a disconnected domain (for example, it may represent grains inside the porous medium). 
    We are interested  in the fluid flow and scalar transport in $\hat{\mathcal{Y}}_f$, neglecting transport within the solid region $\hat{\mathcal{Y}}_s$, that is here represented only through its interface $\hat{\Gamma}$. We do not discuss the upscaling of systems with full conjugate transfer, which has been extensively studied  in a previous work \citep{Municchi2020}.
    If such flow is well described by the Stokes equations (i.e., incompressible, low Reynolds number), it has been shown that the homogenisation procedure leads to the Darcy equation \cite{auriault1995taylor}. Therefore, in this work we will only deal with the scalar transport problem, assuming that the velocity field is prescribed and the upscaled Darcy equation valid.
    We limit here to the case of periodic and non-moving porous media. However, some of these ideas can be extended to stochastic stationary multiscale random media, with Fourier/Bloch \cite{fadili2003stochastic,cherdantsev2019stochastic} or numerical sampling \cite{icardi2016predictivity} approaches, and to slowly varying or quasi-periodic media \cite{van2011homogenisation,brown2011homogenization}. These extensions will be considered in future studies.

\subsection{Reactive transport in porous media}
\label{SS:ADR}
    Let us consider a generic (dimensional) scalar field $\hat{c}\ofxt$ (e.g., solute concentration, temperature) defined in the fluid region $\hat{\Omega}_f$ given by the $\hat{\mathcal{Y}}_f$ of each cell  that obeys the advection-diffusion equation:
    \begin{equation}
    \label{eq::ADEdimensional}
    \begin{dcases}
    \dfrac{\partial \hat{c}}{\partial \hat{t}} +  \hat{\grad}\cdot\left(\hat{\bv}\hat{c}\right) = \hat{\grad} \cdot \left( \diff \hat{\grad} \hat{c} \right)  & \hat{\bx} \; \in \; \hat{\Omega}_f \\
    \diff\grad \hat{c}\cdot \bm{n}=-\kappa\hat{c} + \hat{g}  \; & \hat{\bx} \; \in \; \hat{\Gamma}
    \end{dcases} 
    \end{equation}
    where $\diff$  is the molecular diffusion coefficient, $\hat{\bv}(\hat{\bm{x}})$ a solenoidal (divergence-free) velocity field and $\hat{g}(\hat{\bm{x}})$ is a known forcing (source) term at the boundary.
    The mixed boundary condition on the solid surface $\hat{\Gamma}$  is representing a linear reaction (or surface deposition). For the limit  of the reaction rate $\kappa\to\infty$, it is equivalent to a homogeneous Dirichlet condition.

    Equation~\ref{eq::ADEdimensional} should be provided with a proper set of \emph{external} boundary conditions (e.g., inlet/outlet). Since such boundary conditions are specific for each problem and should not concern the homogenisation procedure, we will not explicitly state them. However, it is important to notice that homogenisation fails when the external boundary conditions play a significant role at the scale $\ell$, since the model with periodic unit cells $\hat{\mathcal{Y}}$ would fail to be a realistic and mathematically consistent description of the porous medium as outlined in \citet{auriault1995taylor}.

We now introduce a set of dimensionless quantities:
    \begin{equation}
    \bm{x}=\frac{\bm{\hat{x}}}{L}; \;\;
    \bv=\frac{\hat{\bv}}{U} \;\;
    c = \frac{\hat{c}} {c_0} ,\;\;
    \hat{t}=\frac{L^2}{\mathcal{D}}t ,\;\;
    g = \frac{\hat{g}}{\kappa c_0} \,,
    \end{equation}
    where $U$ is the system characteristic velocity and $c_0$ is the characteristic value for particle concentration,
    When the dimensionless quantities are substituted into Eq.~\ref{eq::ADEdimensional}, two dimensionless numbers arise:
    \begin{itemize}
        \item The P\'eclet number:
        \begin{equation}
            \label{eq::peclet}
            \PE = \frac{U\ell}{\mathcal{D}}
        \end{equation}
        representing the ratio between inertial and diffusion time scales at the microscale.
        \item The second Damk\"ohler number:
        \begin{equation}
            \label{eq:::DA}
             \DAII = \frac{\kappa \ell}{\mathcal{D}}
        \end{equation}
        which gives the ratio between the reaction and diffusion time scales at the microscale.
    \end{itemize}
    
    Thus, recalling that $\eps = \ell/L$, we can write the advection-diffusion problem in dimensionless form:
    \begin{equation}\label{eq::ADEdimensionless}
    \begin{dcases}
    \dfrac{\partial c}{\partial t} + \div \left( \eps^{-1}\PE \bv c - \grad c \right) = 0 & \bx \; \in \; \Omega \\
    \grad c\cdot \bm{n}=\eps^{-1} \DAII\left( g - c \right) & \bx \; \in \; \Gamma
    \end{dcases}
    \end{equation}

\subsection{The failure of two-scale asymptotics}
\label{SS:fail}
    We show now how the standard homogenisation approach can fail for the simple advection-diffusion problem~\ref{eq::ADEdimensionless}.
    Following the standard two-scale asymptotic homogenisation method \cite{allaire1989homogenization}, we introduce a "fast" spatial variable
    $
    \by=\frac{\bo{x}}{\eps}
    $
     that induces the following differentiation chain rule:
    \begin{equation}
    \label{eq::derivativeExpansionSpatial}
     \grad c = \nx c + \eps^{-1} \ny c
    \end{equation}

    The concentration field $c$ is then represented as $c\ofxt=c\ofxyt$ and expanded into an asymptotic series of $\eps$
    \begin{equation}\label{eq::Cexpansion}
    c\ofxyt=\sum_{m=0}^\infty \eps^mc_m\ofxyt 
    \end{equation}
      
    Furthermore, we consider the case of equally important advection and diffusion at large scales $\eps^{-1} \PE=\PE_L=\mathcal{O}(1)$ without any loss of generality.
    Substituting into Eq.~\ref{eq::ADEdimensionless} and collecting terms of the same order $\eps$ we obtain:
     \begin{align}\label{eq::expandedDispersionstep3}
    &\eps^{-2} \bigg\{\ny \cdot (\PE \bv c_0 - \ny c_0)\bigg\}+& \notag \\
    &\eps^{-1}\bigg\{ \nx \cdot ( \PE \bv c_0 - \ny c_0 ) -\ny \cdot [ (\nx c_0+\ny c_1- \PE\bv c_1)  ] +& \notag \\
    &\eps^0 \bigg\{ \frac{\partial c_0}{\partial t}  - \nx \cdot [  (\nx c_0+\ny c_1)] - \ny \cdot [(\nx c_1+\ny c_2)] + \notag\\
    &\hspace{30pt}+  \PE \nx \cdot( \bv  c_1)\bigg\} 
     = \mathcal{O}(\eps)& 
    \end{align}    
     We apply the same expansion to the boundary condition in $\Gamma$ for the homogeneous case $g=0$:
    \begin{equation}\label{eq::expandedBC}
    (\nx+\eps^{-1}\ny) (c_0+\eps^1c_1+\eps^2c_2)\cdot \bm{n}=
    \eps^{-1} \DAII(g-c_0-\eps^1c_1-\eps^2c_2) + \mathcal{O}(\eps^3) \quad  \by  \in  \Gamma 
    \end{equation}

        Collecting terms of order $\mathcal{O}(\eps^{-2})$ in Eq.~\ref{eq::expandedDispersionstep3} and \ref{eq::expandedBC}, we obtain an equation for the leading order $c_0$:
        \begin{equation}\label{eq::Epsilon-2Terms}
        \begin{cases}
        \ny \cdot (\PE \bv c_0 - \ny c_0) = 0 & \by \; \in \; \mathcal{Y}_f \\
        \ny c_0\cdot\bm{n}= \DAII (g-c_0)  & \by \; \in \; \Gamma 
        \end{cases}
        \end{equation}
        
        Generally, formal two-scale asymptotics can be applied to obtain a macroscale equation (independent from the microscale), if the equation for the lowest order is solved by a function independent of $\by$, i.e., $c_0 = c_0 \ofxt$. This allows the subsequent order equations to simplify and  to finally perform volume averaging to obtain a well defined effective concentration.
        
However, the existence of a constant (in $\by$) solution of Eq.~\ref{eq::Epsilon-2Terms} generally relies on a few conditions that are often overlooked:
\begin{itemize}
\item  A solenoidal velocity field ,i.e., $\ny \cdot \bv =0$, While this is trivially satisfied for solutes and very small particles, we have seen in Section~\ref{SS:smoluchowski} that this is not necessarily true for inertial and colloidal particles, which can instead show regions of accumulation.
\item A slow reaction, such as $\DAII=\mathcal{O}(\eps)$ or $\DAII=0$ (impermeable non reactive walls).
\item Null normal velocity at the wall. While this is clearly verified for the fluid, it might not be valid for particles and colloids.
\item Homogeneity of the reaction terms, i.e., $g=0$. This is not always verified. For example, for non-linear surface reactions, a linearisation will generally involve an inhomogeneous component.
\end{itemize}
The last three assumptions result in a simple homogeneous Neumann boundary condition for the leading order equation, compatible with a constant solutions $c_0$. This is clearly not verified whenever there is a dominant reaction or when there is non-zero velocity at the wall (e.g., membranes, filtration, inertial particles).

Without relying on these assumptions, Eq.~\ref{eq::Epsilon-2Terms} might be non-solvable or solved
 only by the trivial solution $c_0 = 0$, when $\DAII>0$, thus removing the leading order term from the asymptotic expansion. Even when $\DAII=0$, due to the non-solenoidal velocity field or compenetrating velocity at the wall, a non-constant solution is found for $c_0$.
This failure of the standard asymptotic expansion is essentially triggered by a lack of separation of scales, i.e. the leading order of the asymptotic expansion has significant micro-scale variations. Therefore, the problem cannot be homogenised by means of this formal expansion. In the next section we will focus on the case of dominant reaction, following the alternative approach proposed by \citet{ALLAIRE2007523}.

\section{Upscaling transport with dominant heterogeneous reactions}
\label{S:homogeneous}
Let us consider here the case of solute transport (e.g., advection-diffusion equation with divergence-free velocity field, null at the walls) with reactive boundary conditions. In order to overcome the problem highlighted above and find an upscaled equation, we apply a decomposition method by which the scalar field is rewritten as a product of terms that account for exchange processes at different scales.
In order to shorten our notation, we introduce the steady advection-diffusion operator
\begin{equation}
    \label{eq::L}
    \mathcal{L} = \div \left[ \eps^{-1} \PE\bv \ofx\ -  \grad \right],
\end{equation}
that allows us to rewrite equation~\eqref{eq::ADEdimensionless} in a more succinct way:
\begin{equation}
    \label{eq::Caucy_adr}
    \begin{dcases}
   \left( \ddt{} + \mathcal{L} \right) c\ofxt = 0,\,& \quad \forall \bm{x} \in \Omega
   \\
     \grad c\ofxt \cdot \bm{n} = \eps^{-1} \DAII \left(g\ofx-c\ofxt\right) , & \quad  \forall \bm{x}   \in \Gamma
    \end{dcases}
\end{equation} 
As already stated, we do not consider here other ``external" boundary conditions, as if we were considering a ``bulk" region where the dynamics is completely determined by the boundary conditions in $\Gamma$.

In this section we will focus on the approach of \citet{allaire2010two} and \citet{Mauri_1991}, which assumes homogeneous boundary conditions, e.g., $g=0$. In Section~\ref{S:extensions} we will consider general inhomogeneous Robin boundary conditions, as well as non-solenoidal velocities.

\subsection{Spectral decomposition}
\label{SS:decomp}
    By using the geometrical periodicity of the system, we can separate $c\ofxt$ into a function $\phi\ofxt$ periodic in $\mathcal{Y}$ and another function $\omega\ofxt$, which is not necessarily periodic:
    \begin{equation}
        \label{eq::decomposition}
        c \ofxt =\phi\ofx  \omega \ofxt   
    \end{equation}
    The aim of decomposition~\ref{eq::decomposition} is to separate the part of $c \ofxt$, that dominates at the microscale, and drives the mass exchange between $\Omega$ and $\Gamma$, from the part that is responsible for the transport over long distances.
    We therefore choose $\phi$ as the principal (i.e. first positive) eigenfunction of $\mathcal{L}$ in $\mathcal{Y}$, which satisfies:
    \begin{equation}
        \label{eq::phi}
        \begin{dcases}
            \mathcal{L}\phi\ofx =\eps^{-2}\lambda \phi\ofx, & \quad \forall \bm{x} \in \mathcal{Y}_f
       \\
        \grad \phi\ofx \cdot \bm{n}= -\eps^{-1} \DAII \phi\ofx, & \quad  \forall \bm{x} \in \Gamma
        \end{dcases}
    \end{equation}
with $\phi>0 \forall \bx$, and $\lambda$ being the associated principal eigenvalue\footnote{The existence of a real non-negative principal eigenpair is guaranteed by the compactness and positivity of the operator $\mathcal{L}$ even when, like in this case, it is not self-adjoint.
} that, in this case, represents an effective (bulk) reaction constant to balance the reaction at the boundaries.
    Therefore, $\omega \ofxt$ can be thought as a \textit{scaled field}, which compensates for the arbitrariness of our choice of $\phi$.
    
    Notice that $\lambda$ scales as $\eps^{-2}$. This can be easily demonstrated by recalling that the volumetric effects of the boundary reactions, integrating Eq.~\ref{eq::phi}, leads the the following scaling for $\lambda$:
    \begin{equation}
        \label{eq::lambda_eps}
       \eps^{-1}\DAII \frac{\eps^2\int_{\Gamma} \phi \text{d}^2\by}{\eps^3\int_{\mathcal{Y}_f} \phi \text{d}^3\by} = \eps^{-2}\DAII \frac{\int_{\Gamma} \phi \text{d}^2\by}{\int_{\mathcal{Y}_f} \phi \text{d}^3\by} = \eps^{-2} \lambda\,,
    \end{equation}
    where the $\eps$ multiplying the integrals arise from the adimensionalisation of the surface $\Gamma$ and the volume $\mathcal{Y}_f$. Clearly, since $\phi$ is assumed periodic in $\mathcal{Y}_f$, the proper length scale for scaling Eq.~\ref{eq::phi} would be $\ell$. However, it is useful to keep track of the scaling with respect to $L$ for comparison with Eq.~\ref{eq::Caucy_adr}. Furthermore, it is straightforward to see that $\eps$ can be removed from Eq.~\ref{eq::phi} when the differential operators are decomposed according to Eq.~\ref{eq::derivativeExpansionSpatial}.

    Substituting equation~\eqref{eq::decomposition} into $\mathcal{L}c$ we would obtain an equation for $\omega$. However, it is left to the reader to check that this would not lead to a particularly simple expression for $\omega$ and it would not be enough to obtain a proper macroscopic equation.
%
    Therefore, to recast the problem into a more familiar and tractable operator, we introduce the adjoint equation for the adjoint eigenfunction $\phi^{\dagger}$:
     \begin{equation}
        \label{eq::phi_adj}
        \begin{dcases}
            \mathcal{L}^{\dagger}\phi^{\dagger}\ofx = \eps^{-2}\lambda^{\dagger} \phi^{\dagger}\ofx, & \quad \forall \bm{x} \in \mathcal{Y}_f     \\
        \grad \phi^{\dagger}\ofx \cdot \bm{n} = -\eps^{-1} \DAII \phi^{\dagger}\ofx, & \quad  \forall \bm{x} \in \Gamma      
        \end{dcases}
    \end{equation}  
    where the adjoint operator $\mathcal{L}^{\dagger}$ is defined as\footnote{
    It can be proven \citep{lax2014functional,du2006order} that there exist one single real and positive eigenfunction for both the direct and adjoint problem, and the associated eigenvalue $\lambda^\dagger=\lambda$ is real. Since here the velocity field is divergence free and null at the walls, the adjoint problem is simply the original problem with a different sign in front of the advective term. In Section~\ref{S:extensions} we will show how this can be extended.
}
    $$\mathcal{L}^{\dagger}=- \left[  \eps^{-1} \PE\bv\ofx \cdot \grad  +  \div \grad \right]\,.$$
The adjoint problem is related to the infinitesimal generator of the underlying stochastic process describing the system, and it is the operator that evolves observables of the system \cite{Pavliotis2008}. This means its solution is a good measure to weight the concentration field. In our case, we can in fact multiply the equation by the adjoint eigenfunction $\phi^{\dagger}$ to obtain a simplified equation for $\omega$:
    \begin{align}        
       \phi^{\dagger}  \mathcal{L} \left(\phi \omega \right) 
    &=\phi^{\dagger}  \div \left( \eps^{-1} \PE\bv\phi \omega -  \phi \grad  \omega  -  \omega \grad \phi   \right) \nonumber \\
        &= \div \left( \eps^{-1} \PE\bv\phi\phi^{\dagger}\omega - \omega \phi^{\dagger}  \phi \grad  \omega  - \phi^{\dagger}  \grad \phi\omega  \right)  \nonumber \\
        & \quad+ \left(  \phi \grad  \omega  +  \omega \grad \phi   \right)\cdot \grad \phi^{\dagger} - \omega\phi \div \left( \eps^{-1} \PE\bv \phi^{\dagger}\right)
                   \nonumber \\
        &= \div \left( \eps^{-1} \PE\bv\phi\phi^{\dagger}\omega - \omega \phi^{\dagger}  \phi \grad  \omega  - \phi^{\dagger}  \grad \phi\omega  + \grad\phi^{\dagger}   \phi\omega\right) \nonumber \\ 
        &\quad- \omega\phi \div \left( \eps^{-1} \PE\bv \phi^{\dagger} + \grad \phi^{\dagger} \right) \nonumber \\
        &= \div \left( \bv^{\star} \omega-   \beta\grad \omega \right) +\eps^{-2} \lambda\beta\omega
    \end{align}
    where, in the last step, we have used the definition of $\mathcal{L}^{\star}$, and we have introduced a dimensionless coefficient $\beta \ofx = \phi \phi^{\dagger}$ and a new velocity field:
    \begin{equation}
        \label{eq::vstar}
        \bv^{\star} =  \beta\left( \bv  - \frac{\grad \phi}{\PE \phi} + \frac{\grad \phi^{\dagger} }{\PE \phi^{\dagger}}\right) 
    \end{equation}
%
Notice that $\grad \phi \sim \mathcal{O}\left(\eps^{-1} \right)$ since $\phi$ (and $\phi^{\dagger}$) are periodic in $\mathcal{Y}_f$.

    We can now define a new advection-diffusion operator:
    \begin{equation}
        \mathcal{L}^{\star} = \div \left[ \eps^{-1} \PE \bv^{\star} \ofx - \beta \ofx  \grad \right]
    \end{equation}
    and an equation for $\omega$ that reads:
    \begin{equation}
        \label{eq::omega_eq}
        \left[\beta\left(\ddt{} + \frac{\lambda}{\eps^{2}}\right) + \mathcal{L}^{\star}\right]\omega = 0 
    \end{equation}
    
    The modified operator $\mathcal{L}^{\star}$ possesses several good properties. In particular, the new velocity field $\bv^{\star}$ is divergence-free:
        \begin{align}
            \div \bv^{\star} &= \div \left( \bv\phi\phi^{\dagger} + \grad\phi^{\dagger} D  \phi -  \phi^{\dagger} D \grad \phi \right) \nonumber \\
            &= \phi\div \left( \bv \phi^{\dagger}  + D\grad\phi^{\dagger}\right) + \phi^{\dagger}\div \left( \bv \phi  - D\grad\phi\right) \nonumber \\
            &= -\phi\phi^{\dagger}\frac{\lambda}{\eps^2} + \phi\phi^{\dagger}\frac{\lambda}{\eps^2} = 0
        \end{align}        

If we apply the same decomposition to the boundary condition in $\Gamma$, it can be easily verified that, due to our choice of $\phi$, the Robin (reactive) boundary condition is absorbed by $\phi$, leaving a homogeneous Neumann  (non-reactive) BC for $\omega$:
    \begin{align}
        \label{eq::w_BC}
        \grad\omega \cdot \bm{n}  = 0
    \end{align}

    Therefore, the spectral decomposition has the effect of moving the reaction term from the boundary to the bulk. It is easy to see now that equation~\eqref{eq::omega_eq} can be further simplified by means of one additional transformation which takes into account the fast time scale related to $\lambda$:
    \begin{equation}
        \label{eq::dec_time}
        \omega \ofxt = e^{-\eps^{-2}\lambda t} w \ofxt \to  c \ofxt =\phi\ofx  e^{-\eps^{-2}\lambda t} w \ofxt  \,,
    \end{equation}
    where the exponential takes into account the fast change in $c$ due to the reaction.
     We can finally deduce a simple transport equation for the reduced concentration $w$:
    \begin{equation}
        \label{eq::w_eq}
        \begin{dcases}
         \left(\beta\ddt{} + \mathcal{L}^{\star} \right) w = 0 ,& \quad \forall \bm{x} \in \Omega
        \\
        \grad w \cdot \bm{n} = 0, & \quad \forall \bm{x} \in \Gamma
        \end{dcases}
    \end{equation}
Through this simple decomposition, we have obtained an equation for the reduced concentration with non-reactive (zero flux) on $\Gamma$. This can be now upscaled using standard two-scale asymptotics.

\subsection{Two-scales asymptotics with drift}

Decomposition~\eqref{eq::dec_time} does not rely on any scale separation or two-scale expansion. However, due to the periodicity, we considered auxiliary problems for the eigenfunction $\phi$ and the adjoint $\phi^{\dagger}$ in the unit cell $\mathcal{Y}_f$. They can therefore only depend on the small scale spatial variable $\bm{y}$, while $w$ can vary on both scales. It is now convenient to express the solution using the standard two-scale formalism:
\begin{equation}
    \label{eq::c_multiscale}
    c \left( \bm{x} ,\bm{y} , t\right) = e^{-\eps^{-2}\lambda t} \phi \left( \bm{y}\right) w \left( \bm{x} ,\bm{y} , t\right)
\end{equation}

The two-scale expansion of $c$ is fully determined by expansion of $w$. Since the governing equation for $w$ no longer has a reactive boundary condition,  the limitation of standard asymptotics addressed in Section~\ref{S:intro}, can be avoided.
We thus expand $w$ into an asymptotic series of $\varepsilon$ introducing a drift in $w$ that accounts for the fast advective transport at the microscale:
\begin{align}\label{eq::wexpansion}
&w(\bm{x}-\eps^{-1}  \PE\bV^{\star}t ,\bm{y},t)=\sum_{m=0}^\infty \eps^mw_m(\bo{x} - \eps^{-1}  \PE\bV^{\star}t,\by,t) &
\end{align}
where the velocity $\bV^{\star}$ is obtained during the homogenisation procedure as shown in \cite{allaire2010two}.
The same expansion we derived in Section~\ref{S:intro} is applied on Eq.~\ref{eq::w_eq} obtaining:
\begin{align}\label{eq::expanded_w}
&\eps^{-2} \bigg\{\ny \cdot (\PE \bv^{\star} w_0 - \beta\ny w_0)\bigg\}+& \notag \\
& \eps^{-1}\bigg\{   \PE\left[-\beta\bV^{\star}\cdot \nx w_0   +\nx \cdot (\bv^{\star}  w_0)+\ny \cdot (\bv^{\star}  w_1)\right] \-& \notag \\
&\hspace{15pt}- \nx \cdot (\beta \ny w_0 ) -\ny \cdot [ \beta(\nx w_0+\ny w_1) ] \bigg\}+ & \notag \\
&+\eps^0 \bigg\{ \beta \ddt{w_0} - \beta \PE \bV^{\star} \cdot \nx w_1 + \PE\nx \cdot (\bv^{\star} w_1)+ \PE\ny \cdot (\bv^{\star} w_2)   -&  \notag \\
&\hspace{15pt}- \nx \cdot [ \beta  (\nx w_0+ \ny w_1)] - \ny \cdot [\beta(\nx w_1+\ny w_2)]\bigg\}= \mathcal{O}(\eps)
\end{align}
The internal boundary condition is also expanded:
\begin{align}\label{eq::expandedwBC}
&(\nabla_x+\eps^{-1}\nabla_y) (w_0+\eps^1w_1+\eps^2w_2)\cdot \bm{n}=0
\end{align}
which corresponds to a homogeneous Neumann boundary condition for all the terms in the series.

\subsubsection*{Terms of order $\mathcal{O}(\eps^{-2})$}
    We can now collect the leading order terms of equation~\eqref{eq::expanded_w} and \eqref{eq::expandedBC}, which correspond to terms of order $\mathcal{O}(\eps^{-2})$
    \begin{align}\label{eq::wEpsilon-2Terms}
    \begin{cases}
    \ny \cdot (\PE \bv^{\star} w_0 - \beta\ny w_0) =0  & \by \; \in \; \Omega \\
    \ny w_0\cdot\bm{n}=0  & \by \; \in \; \Gamma 
    \end{cases}
    \end{align}
    
    This problem allows the trivial solution $w_0 = w_0 \left(\bm{x},t\right)$ since the homogeneous Neumann boundary condition ensures that no terms from "fast" scales appear in equations~\eqref{eq::wEpsilon-2Terms}.

\subsubsection*{Terms of order $\mathcal{O}(\eps^{-1})$}
At the order $\mathcal{O}(\eps^{-1})$, we obtain the following partial differential equation:
\begin{align}\label{eq::wEpsilon-1Terms}
\begin{cases}
  \PE\beta\bV^{\star}\cdot \nx w_0   + \PE\nx \cdot (\bv^{\star}  w_0)
 + & \\
 \quad + \ny \cdot [ \PE\bv^{\star}  w_1 -  \beta(\nx w_0+\ny w_1) ]  & \by \; \in \; \Omega \\
\left( \ny w_1 + \nx w_0 \right)\cdot\bm{n}=0  & \by \; \in \; \Gamma 
\end{cases}
\end{align}

This equation can be integrated over the unit cell to give:
\begin{equation}
    \label{eq::wEpsilon-1integrated}
    \left[\int \limits_{\mathcal{Y}_f} \left(\beta \bV^{\star} -  \bv^{\star}\right) \text{d}\bm{y}\right] \cdot \nx w_0
    = 0
\end{equation}

Without loss of generality\footnote{Since both the direct and adjoint eigenfunctions can be defined up to a constant.}, we can take $\beta$ such that:
\begin{equation}
    \label{eq::norm_beta}
    \int \limits_{\mathcal{Y}_f} \beta \text{d}\bm{y} = \int \limits_{\mathcal{Y}_f}  \text{d}\bm{y} = \epsilon
\end{equation}

Where $\epsilon \in \left[0,1 \right]$ is the porosity of the porous medium, and is a consequence of our choice of dimensionless variables. In fact, the maximum volume that the fluid can occupy inside the unit cell is $\ell^3$, which corresponds to $\epsilon = 1$.

It is possible to employ the porosity to define an averaging operator that will be employed in the following upscaling:

\begin{equation}
    \label{eq::Favre_avg}
    \left\langle \star \right\rangle = \frac{1}{\epsilon} \int \limits_{\mathcal{Y}_f} \star \; \text{d}\bm{y}
\end{equation}

This corresponds to the Favre averaging operator \cite{A.Favre1965} often employed in the description of compressible turbulent flows and multiphase flows \cite{Radl2018}.
Carrying on the integration in Eq.~\ref{eq::wEpsilon-1integrated}  results in an expression for $V^{\star}$:
\begin{align}
      \label{eq::ddta_w0}
      \bV^{\star} &= \frac{1}{\epsilon} \int \limits_{\mathcal{Y}_f}  \bv^{\star}(\bm{y})\text{d}\bm{y}  = \langle\bv^{\star} \rangle \nonumber\\
\end{align}


Combining \eqref{eq::ddta_w0} and equation~\eqref{eq::wEpsilon-1Terms}, we can express the solution $w_1$ as:
\begin{equation}
    \label{eq::first_ord_corrector}
    w_1 = \boldsymbol{\chi} \ofy \cdot \nx w_0 + f_1 \left( \bm{x} \right)
\end{equation}
where $\boldsymbol{\chi}$ is called first order corrector and $f_1$ is an arbitrary function of the macroscopic coordinate.
Substituting equation~\eqref{eq::first_ord_corrector} into equations~\eqref{eq::wEpsilon-1Terms} we obtain:
\begin{equation}
    \label{eq::first_order_corrector_eq}
    \begin{dcases}
    -\ny\cdot \left[ \beta \left( \mathbf{I} + \ny\boldsymbol{\chi} \right) \right] +  \PE\bv^{\star} \cdot \left( \mathbf{I} + \ny \boldsymbol{\chi}\right)=  \PE\beta\bV^{\star}  & \by \; \in \; \mathcal{Y}_f \\
    \left(\mathbf{I} + \ny \boldsymbol{\chi} \right) \cdot \bm{n} = 0& \by \; \in \; \Gamma
    \end{dcases}
\end{equation}

Equation~\eqref{eq::first_order_corrector_eq} is often referred to as the first order cell corrector, closure, or cell problem, since $\boldsymbol{\chi}$ appears in the expressions for the effective diffusion coefficient.
It is important to notice that that the first order corrector only depends on $\by$ and therefore, its governing equation can be defined in the unit cell $\mathcal{Y}$ alone.

\subsection{Upscaled equation}

Finally, we collect the terms of order $\mathcal{O}(1)$ from equation~\eqref{eq::expanded_w}:
\begin{align}\label{eq::wEpsilon-0Terms}
\begin{dcases}
\beta \ddt{w_0} + \beta \PE (\bv^{+} - \bV^{\star}) \cdot \nx w_1 +  + \PE\ny \cdot (\bv^{\star} w_2)  -&   \\
\hspace{15pt}- \nx \cdot [ \beta  (\nx w_0+ \ny w_1)]  - \ny \cdot [\beta(\nx w_1 + \ny w_2)] & \by \; \in \; \Omega \\
\left( \nx w_1 + \ny w_2 \right)\cdot\bm{n}=0  & \by \; \in \; \Gamma
\end{dcases}
\end{align}

equation~\eqref{eq::wEpsilon-0Terms} can be integrated over $\mathcal{Y}_f$ to obtain:
\begin{align}
    \label{eq::macroEq1}
       \langle \beta \rangle \ddt{ w_0 }  - \nx \cdot \left[\langle \beta\left(\mathbf{I} + \ny \boldsymbol{\chi} +  \PE\bV^{\star}\boldsymbol\chi - \PE\bv^{+}\boldsymbol\chi  \right)\rangle \cdot \nx w_0  \right] =0
\end{align}

Furthermore, we can express an effective diffusivity tensor $\bm{D}^{\text{eff}}$ as:
\begin{equation}
    \label{eq::Deff_1}
    \bm{D}^{\text{eff}} = \langle \beta\left(\mathbf{I} + \ny \boldsymbol{\chi} +  \PE\bV^{\star}\boldsymbol\chi -  \PE\bv^{+}\boldsymbol\chi  \right)\rangle
\end{equation}

The full form of the effective diffusion tensor is irrelevant since its effect on the homogenised solution is defined by its product with the Hessian matrix $\nx^2 w_0$, so that only its symmetric part contributes in the equation. It is thus more convenient to express it in symmetric form. It can be shown \citep{auriault1991} that the symmetric part of $\bm{D}^{\text{eff}}$ can be expressed as:
\begin{equation}
    \label{eq::Deff_2}
    \bm{D}^{\text{eff}} = \langle \beta\left(\mathbf{I} + \ny \boldsymbol{\chi}\right) \cdot \left(\mathbf{I} + \ny \boldsymbol{\chi} \right)^{T} \rangle
\end{equation}

This implies that $\bm{D}^{\text{eff}}$ does not depend on the microscopic velocity explicitly, but through $\boldsymbol{\chi}$. 

Finally, to include the effect of the exponential term in the decomposition of $w$, we transform back the shifted spatial coordinate, and we let $\eps \to 0$ to obtain the homogenised equation for $\omega\ofxt = \omega_0 \ofxt + \mathcal{O}(\eps)$:
 \begin{equation}
    \label{eq::macroscopic_equation}
    \ddt{\omega}  + \nx  \cdot \left(\PE_L\bV^{\star}\omega - \bm{D}^{\text{eff}} \cdot \nx \omega \right) =- \frac{\lambda}{\eps^2} \omega
\end{equation}
where we  have introduced the macroscopic P\'eclet number $\PE_L = UL/\mathcal{D} = \eps^{-1} \PE$, which is evaluated at the macroscale. Similarly, the reactive term is still of order $\mathcal{O}(\eps^{-2})$. 
It should be noted that equation~\eqref{eq::macroscopic_equation} is not immediately equivalent to the equation for the average concentration $\langle c \rangle$. In fact:
\begin{equation}
    \label{eq::ave_c}
    \langle c \rangle = \langle \phi \omega \rangle = \frac{1}{\epsilon}\int \limits_{\mathcal{Y}_f} \phi \ofy \omega (\bm{x},\bm{y},t) \text{d}\bm{y}
\end{equation}

However, employing expansion \eqref{eq::wexpansion} and letting $\eps \to 0$, at the leading order we have:
\begin{equation}
    \langle c \rangle =  \frac{1}{\epsilon}\int \limits_{\mathcal{Y}_f} \phi \ofy \omega \ofxt \text{d}\bm{y} = \langle\phi \rangle \omega_0 \ofxt + \mathcal{O}(\eps)
\end{equation}

Finally, taking a $\phi$ normalised over the fluid volume\footnote{There is now a unique choice of $\phi$ and $\phi^\dagger$ such that $ \langle\phi \rangle=1$ and $ \langle\phi\phi^\dagger \rangle=1$.}, we can write the asymptotic limit:
\begin{equation}
    \label{eq::cave_1}
    \langle c \rangle \sim \omega_0 \ofxt, \quad \eps \to 0.
\end{equation}
and therefore rewrite a macroscopic equation for the first order approximation\footnote{ The fact that we took the limit $\eps \to 0$ is equivalent to state that we employed a local perturbation analysis. This has the consequence that Eq.~\ref{eq::cave_equation} is expected to be valid only in the case in which the microscopic and macroscopic scales are totally separated and $\ell \ll L$.} of $\langle c \rangle$:
\begin{equation}
    \label{eq::cave_equation}
    \ddt{ \langle c \rangle}  + \grad  \cdot \left(\PE_L\bV^{\star} \langle c \rangle- \bm{D}^{\text{eff}} \cdot \grad  \langle c \rangle \right) =- \frac{\lambda}{\eps^2}  \langle c \rangle
\end{equation}

Alternatively, one can retain terms up to $\mathcal{O}(\eps^2)$ and merely consider equation~\eqref{eq::macroscopic_equation} as the equation for $w_0$, i.e. the leading order term.
The average concentration can be corrected a-posteriori by introducing a first order term:
\begin{equation}
    \label{eq::cave_2}
    \langle c \rangle = \omega_0 + \eps\langle \phi \boldsymbol{\chi} \rangle \cdot \grad  \omega_0 + \mathcal{O}\left( \eps^2\right), \quad \quad \eps < 1
\end{equation}

Equation~\ref{eq::cave_2} is often referred as the ``first corrector equation" \cite{Davit2013} but a direct solution would involve a third order macroscopic equation that is therefore not easily solvable.

Notice that $\langle c \rangle$ is a smooth field at the microscale, since the homogenisation procedure completely removed the high wave number modes relative to the scales smaller then $\ell$. Thus, Eq.~\ref{eq::cave_equation} can be rescaled to the lenght $\ell$ by setting $\eps = 1$. This is similar to what is commonly done in the method known as Renormalisation Group \citep{Wilson1974,Forster1977}, where the wave number is rescaled to compensate for the loss of degrees of freedom after removing the high wave number components from the microscopic equations. This freedom reflects the obvious fact that the final model cannot depend on the choice of the reference length scale. In other words, we can rescale Eq.~\ref{eq::cave_equation} over $\ell$ rather than $L$ without any loss of generality, resulting in:

\begin{equation}
    \label{eq::cave_eq_rescaled}
    \ddt{ \langle c \rangle}  + \grad  \cdot \left(\PE\bV^{\star} \langle c \rangle- \bm{D}^{\text{eff}} \cdot \grad  \langle c \rangle \right) =- \lambda  \langle c \rangle    
\end{equation}

Reverting Eq.~\ref{eq::cave_equation} to its dimensional form results in the following partial differential equation for the dimensional field $\langle c \rangle$:

\begin{equation}
    \label{eq::cave_eq_dimensional}
    \frac{\text{d} \langle \hat{c} \rangle }{\text{d} \hat{t}} + \hat{\grad} \cdot \left(\hat{V}^{\star}   \langle \hat{c} \rangle - \hat{\bm{D}^{\text{eff}}} \hat{\grad}\langle \hat{c} \rangle \right) = \hat{\lambda} \langle \hat{c} \rangle
\end{equation}

Where the new dimensional parameters are defined as:
\begin{equation}
    \label{eq::dim_param}
    \hat{V}^{\star} =V^{\star} U, \quad \hat{\bm{D}^{\text{eff}}} = \bm{D}^{\text{eff}} \mathcal{D}, \quad \hat{\lambda} = \frac{\lambda \mathcal{D}}{\ell^2} 
\end{equation}

Notice that $\hat{\lambda}$ is obtained from the dimensionless eigenvalue by means of the length scale $\ell$. This is due to the order $\eps^{-2}$, which reflects the fact that the cell problem is normally made dimensionless with the length scale $\ell \neq L$.

\section{Inhomogeneous boundary conditions and electrostatic forces}
\label{S:extensions}

\subsection{Extension to inhomogeneous boundary conditions}
The analysis presented in Section~\ref{S:homogeneous} has been first introduced by \citet{Allaire1992} and \citet{Mauri_1991} and works well for the case of homogeneous Robin boundary conditions (i.e., when $g = 0$ or when $g$ is a linear function of $c$ and $\grad c \cdot \bm{n} $  ), but it loses its generality when applied to the general inhomogeneous case defined by equations~\ref{eq::Caucy_adr}.

In fact, applying the decomposition $c = \omega \phi$  results in the following boundary condition:
\begin{equation}
    \label{eq::Robinomegaphi}
    \phi \grad \omega \cdot \bm{n} = \omega \left(-\eps^{-1} \DAII\phi - \grad \phi \cdot \bm{n} \right) + \mathcal{G} \ofx \,,
\end{equation}

where we set $\mathcal{G}=\eps^{-1}\DAII g$ to simplify the subsequent analysis.

Equation~\ref{eq::Robinomegaphi} gives rise to a new inhomogeneous boundary condition for $\omega$:

\begin{equation}
    \label{eq::omegaBC1}
    \grad \omega \cdot \bm{n} =  \phi^{-1} \mathcal{G}\,.
\end{equation}

It is possible to identify three scenarios based on the scaling of $\phi^{-1} \mathcal{G} \ofx$ with $\varepsilon$.

\begin{itemize}
    
\item $\phi^{-1} \mathcal{G} \ofx \approx \mathcal{O}\left(\varepsilon^{-1}\right)$

In most applications,  $\phi^{-1} \mathcal{G}$ is a function of $\eps^{-1} \DAII$ and varies in $\mathcal{Y}_f$. This is the standard case where the boundary condition for $c$ can be written as:
\begin{equation}
    \label{eq::BCc1}
    \grad c \cdot \bm{n} = \eps^{-1} \DAII \left[ c - 
    g\ofx \right] \,,
\end{equation}
where $g\ofx$ is some function defined on $\Gamma$ representing, for example, another concentration field. It is easy to see that applying the method illustrated in Section~\ref{S:homogeneous} to a boundary condition as Eq.~\ref{eq::BCc1}, leads to a hierarchy of equations suffering of the same pathological problem described in Section~\ref{SS:fail} (i.e., the leading order term of the $c$ expansion still depends on the $\bm{y}$ coordinate). 
This issue frustrates the upscaling procedure whenever $\phi^{-1} \mathcal{G}$ varies at the same scale of the microscopic gradients, and thus, when it varies as $\varepsilon^{-2}$. Unfortunately, this is also a natural scaling for fluxes, that are defined per unit area of the microscopic boundary.

\item $\phi^{-1} \mathcal{G} \ofx \approx \mathcal{O}\left(1\right)$

Even when $\phi^{-1} \mathcal{G}$ is defined per unit length, and thus scales as $\varepsilon^{-1}$, the method fails. It is easy to show that in this case it is not possible to perform decomposition~\ref{eq::first_ord_corrector} anymore, and thus $c_1 = c_1 \left(\bm{x},\bm{y}\right)$ is not periodic in $\mathcal{Y}_f$. This makes the corrector problem defined in all $\Omega$, and thus removes any advantage of performing homgenisation in comparison to solving the original problem.

\item $\phi^{-1} \mathcal{G} \ofx \approx \mathcal{O}\left(\eps\right)$

Finally, it is easy to verify that in the simple case where  $\phi^{-1} \mathcal{G}$ scales as $\varepsilon$, the equation for $\left\langle c \right\rangle$ is simply supplemented with an additional source term, resulting in:

\begin{equation}
    \label{eq::cave_Jweak}
    \ddt{ \langle c \rangle}  + \grad  \cdot \left(\PE_L\bV^{\star} \langle c \rangle- \bm{D}^{\text{eff}} \cdot \grad  \langle c \rangle \right) =- \eps^{-2}\lambda  \langle c \rangle - \int_{\Gamma}\phi^{\dagger}  \mathcal{G} \text{d}\by \,.
\end{equation}

This can be considered a case of adjoint homogenisation, where  the source term is weighted with the adjoint function. Notice that this is similar to what arises in nuclear engineering applications of homogenisation, where $\phi^{\dagger}$ is called ``importance function'' \cite{Stacey2007}.

\end{itemize}

Therefore, it can be concluded that the upscaling procedure based on the works of \citet{Allaire1992} and \citet{Mauri_1991} does not extend to the general case of inhomogeneous Robin boundary conditions.
%
%
It is possible to circumvent the aforementioned problems in the case when $\phi^{-1} \mathcal{G}$ is periodic in $\mathcal{Y}$ by introducing an auxiliary problem. The original problem is separated into a homogeneous problem and an inhomogenous auxiliary problem. Thus, the field $c\ofxt$ is linearly decomposed as follows:
\begin{equation}
    \label{eq::decomplin}
    c \ofxt = c_h \ofxt + \psi \ofx \,,
\end{equation}
where $c_h \ofxt$ satisfies a homogeneous problem with homogeneous boundary conditions:
\begin{equation}
    \label{eq::Caucy_adr_ch}
    \begin{dcases}
   \left( \ddt{} + \mathcal{L} \right) c_h \ofxt =  0,\,& \quad \forall \bm{x} \in \Omega
   \\
     \grad c_h \ofxt \cdot \bm{n} = \eps^{-1} \DAII c_h \ofxt , & \quad  \forall \bm{x}   \in \Gamma\,.
    \end{dcases}
\end{equation} 

While $\psi \ofx$ satisfies a homogeneous auxiliary problem with inhomogeneous boundary conditions, under the assumption that  $\phi^{-1} \mathcal{G}\ofx$ is periodic in $\mathcal{Y}_f$:
\begin{equation}
    \label{eq::aux_problem}
    \begin{dcases}
   \mathcal{L}\psi \ofx = 0,\,& \quad \forall \bm{x} \in \mathcal{Y}_f
   \\
     \grad \psi \ofx \cdot \bm{n} = \eps^{-1} \DAII \psi \ofx + \mathcal{G}\ofx, & \quad  \forall \bm{x}   \in \Gamma
    \end{dcases}
\end{equation}

Notice that the solvability condition for Eq.~\ref{eq::aux_problem} requires:
\begin{equation}
    \label{eq::aux_solvability}
    \eps^{-1} \int_{\Gamma} \DAII \psi \text{d}\Gamma = -\int_{\Gamma} \mathcal{G} \text{d}\Gamma \,. 
\end{equation}

Therefore, if the boundary flux is in the form employed in Eq.~\ref{eq::BCc1}, the solvability condition requires that:
\begin{equation}
    \label{eq::aux_solvability2}
    \int_{\Gamma} \psi\ofx \text{d}\Gamma = \int_{\Gamma} g\ofx \text{d}\Gamma   \,. 
\end{equation}

Equation~\ref{eq::aux_problem} represents one additional cell problem that needs to be solved on $\mathcal{Y}_f$ for $\psi$. After this linear decomposition, the procedure is identical to that described in details in Section~\ref{S:homogeneous}, with the difference that now $c$ is substituted by $c_h$. Finally, one obtains an equation for $\left\langle c \right\rangle$ substituting $\omega_0 = \left\langle c \right\rangle - \left\langle \psi \right\rangle$ into Eq.~\ref{eq::cave_equation}:
\begin{equation}
    \label{eq::cave_inh}
    \ddt{ \langle c \rangle}  + \grad  \cdot \left(\PE_L\bV^{\star} \langle c \rangle- \bm{D}^{\text{eff}} \cdot \grad  \langle c \rangle \right) = - \frac{\lambda}{\eps^2} \left( \langle c \rangle - \langle \psi \rangle \right) \,.
\end{equation}

Following this approach, the effects of the inhomogeneous term in the boundary condition results in a constant source term in the macroscopic governing equation, and does not affect the upscaled transport coefficients. This is consistent with the fact that we assumed the inhomogeneous term constant in time, but requires the solution of one additional cell problem.

\subsection{Non-solenoidal velocity fields and colloidal particles}
\label{SS:smoluchowski}
Colloidal particles are subject to a number of forces such as gravity  and electrostatic forces \citep{elimelech2013particle} in addition to advection and Brownian motion.
Under the Smoluchowski approximation (small Stokes number), these can be treated as additional fluxes in the equation:
    \begin{equation}
    \label{eq::colloids}
    \begin{dcases}
    \dfrac{\partial \hat{c}}{\partial \hat{t}} +  \hat{\grad}\cdot\hat{\bm{j}} & \hat{\bx} \; \in \; \hat{\Omega}_f \\
    \hat{\bj}\cdot \bm{n}=-\kappa\hat{c} + \hat{g}\; & \hat{\bx} \; \in \; \hat{\Gamma}
    \end{dcases} 
    \end{equation}
The total flux $\hat{\bj}$ can be expressed as the usual advection--diffusion flux, plus the Smoluchowski contribution:
$$
\hat{\bj}= \hat{\bv} \hat{c} - \gamma^{-1}\hat{\grad}\hat{\Lambda} \hat{c} - \diff\hat{\grad} \hat{c} \,,
$$
and $\Lambda$ is the underlying gravitational or electrostatic potential, and $\gamma$ is a friction coefficient. After recasting this equation in dimensionless form, we obtain:
    \begin{equation}\label{eq::ColloidsAdim}
    \begin{dcases}
    \dfrac{\partial c}{\partial t} + \div \bj = 0 & \bx \; \in \; \Omega \\
    \bj\cdot \bm{n}=\eps^{-1} \DAII\left( g - c \right) & \bx \; \in \; \Gamma
    \end{dcases}
    \end{equation}
and, defining a new (dimensionless) velocity $\bv_\Lambda=-\grad\Lambda$ as the gradient of the potential, and the dimensionless number $\mu=\frac{\Lambda_0}{\gamma U \ell}$ as the ratio between a reference potential difference and the fluid response, we can write:
\begin{equation}
    \label{eq::J_smoluchowski}
    {\bj}= \eps^{-1}\PE{\bv} c + \eps^{-1}\PE\mu\bm{v}_\Lambda c - {\grad} c \,.
\end{equation}
It should be noted that $\bm{v}_{\Lambda}$ does not necessarily vanish at the boundaries. Therefore the additional flux $\eps^{-1}\PE \mu \bm{v}_{\Lambda} c$ generates accumulation or rarefaction near the wall, depending on the direction of this flux. This flux however is homogeneous in $c$ and it can be treated like a reaction term employing Robin conditions. The second challenge of this model is that the new velocity is no longer solenoidal, unless we consider a simple potential like the gravitational one.

When applying the two-scale expansion, the potential is generally a function of both $\bx$ and $\by$ (we assume that the potential does not vary in time), so that $\Lambda = \Lambda (\bx,\by)$.
This induces the following decomposition on $\bv_{\Lambda}$:
\begin{equation}
    \label{eq::vLambda}
    \bv_{\Lambda} = - \left( \nx + \frac{1}{\eps} \ny \right) \Lambda 
    = \bv_{\Lambda,\bx} + \frac{1}{\eps}\bv_{\Lambda,\by} 
\end{equation}
The second term in $\bv_{\Lambda}$ scales exactly as $\bv$ and thus they can be summed.
In the case the potential does not vary at the macroscopic scale as in the case of colloids, the procedure presented so far is still valid without any major change, provided that now the proper adjoint operator is considered. Instead of Eq.~\ref{eq::phi_adj}, in this case, the adjoint operator is:
     \begin{equation}
        \label{eq::phi_adj}
        \begin{dcases}
            \grad\cdot\left(-\eps^{-1}\PE{\bv} \phi^{\dagger} - \eps^{-1}\PE\mu\bm{v}_{\Lambda,\by}\phi^{\dagger} - {\grad} \phi^{\dagger}\right)
            +
            \eps^{-1}\PE\mu\left(\grad\cdot\bm{v}_{\Lambda,\by}\right)\phi^{\dagger}
            =
            \eps^{-2}\lambda^{\dagger} \phi^{\dagger},
            &
            \quad \forall \bm{x} \in \mathcal{Y}_f
            \\
        \grad \phi^{\dagger}\ofx \cdot \bm{n}
        =
        -\eps^{-1} \DAII \phi^{\dagger}\ofx, & \quad  \forall \bm{x} \in \Gamma      \,.
        \end{dcases}
    \end{equation}

Compared to the direct problem, Eq.~\ref{eq::ColloidsAdim}, not only has an advection term with different velocity but it also has different boundary condition that do not consider the velocity at the wall and an extra term proportional to the divergence of the velocity.

In the more general case of a potential that varies at the macroscale $\bx$ too, i.e.,  $\bv_{\Lambda,\bx} \neq 0$, the upscaling procedure fails. While our handling of inhomogenous boundary conditions is not affected (the contribution from $\bv_{\Lambda,\bx}$ is linear in $c$) and the eigenproblems remain the same (because the term is of order $\eps$ in the eigenproblems), one must modify the equation for $w$ at order $\eps^{-1}$ (Eq.~\ref{eq::wEpsilon-1Terms}) as follows:
\begin{align}\label{eq::wEpsilon-1Terms_pot}
    \begin{cases}
      \PE\beta\bV^{\star}\cdot \nx w_0   + \PE\nx \cdot (\bv^{\star} w_0)
     + & \\
     \quad + \ny \cdot [ \PE\bv^{\star}  w_1 + \PE\mu\beta\bv_{\Lambda,\bx} w_0 -  \beta(\nx w_0+\ny w_1) ]  & \by \; \in \; \Omega \\
    \left( \ny w_1 + \nx w_0 - \PE\mu\beta \bv_{\Lambda,\bx} w_0 \right)\cdot\bm{n}=0  & \by \; \in \; \Gamma 
    \end{cases}
\end{align}

Notice that a new term proportional to $w_0$ arises which does not allow to define a corrector field $\boldsymbol{\chi}\ofy$. The corrector would have, in fact, an explicit dependence on the macroscopic spatial variable $\bx$. Therefore, in the remainder of this work we will always assume that the potential varies at the microscopic scale only and that $\bv_{\Lambda} = \bv_{\Lambda,\by}$ scales as $\eps^{-1}$.

\section{Numerical implementation of the upscaling method}
\label{S:num}
The upscaling procedure explained above is implemented within the C++ opensource finite volume library OpenFOAM\textsuperscript{\textregistered} \citep{Foundation2014} to solve the closure problems in general geometries. We motivate our choice of OpenFOAM \textsuperscript{\textregistered} over other libraries with its wide diffusion  both in the academic and industrial communities, and with the wide range of classes already available in the library and structured in an consistent object-oriented programming approach.   

%

\begin{figure}[htbp]
    \centering
    \includegraphics{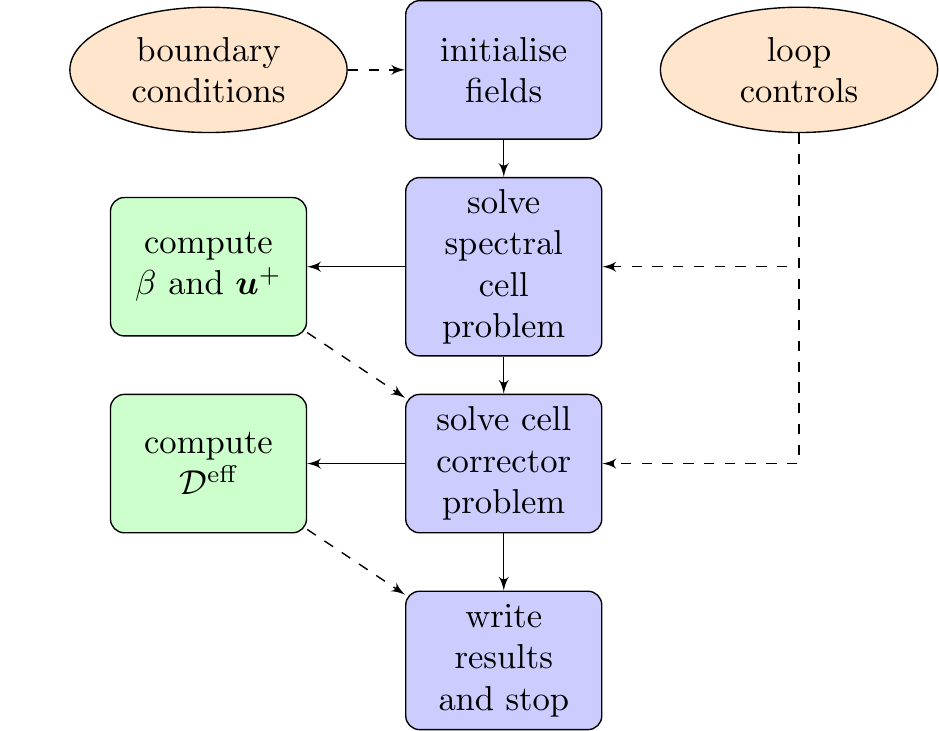}
    \caption{Overview of the numerical procedure. Orange clouds indicate that data is read from OpenFOAM\textsuperscript{\textregistered} dictionaries, blue boxes indicate operations and green boxes indicate the computation of quantities relevant to subsequent operations. Dashed lines represent flow of information.}
    \label{fig::flowchart_tot}
\end{figure}

Figure~\ref{fig::flowchart_tot} illustrates the overall algorithms, which consists of two main sequential operations: first the solution of the spectral cell problem for the direct and adjoint equations, then solving the cell corrector problem for the first order corrector $\boldsymbol{\chi}$.
As input, the algorithm requires an appropriate velocity field which can be obtained from native OpenFOAM\textsuperscript{\textregistered} solvers such as \texttt{simpleFoam}. 

\subsection{Power method for the spectral problem}
 
 Solving the spectral problems poses an additional complication with respect to standard power methods as the direct and adjoint problems are coupled through $\lambda$ and $\lambda^\dagger$ respectively, that should be equal. 
 We propose an iterative segregated algorithm where the convergence of $\lambda$ is achieved through residual control. At each iteration $n$, the values of eigenfunctions and eigenvalues at iteration $n+1$ are calculated following a series of steps:
 
 \begin{itemize}
     \item[1.] Compute $\phi^{n+1}$ and $\phi^{\dagger,n+1}$ from:
     \begin{equation}
        \phi^{n+1} = \mathcal{L}^{-1} \left(\lambda^n \phi^n\right) ,
     \end{equation}
     \begin{equation}
         \phi^{\dagger,n+1}=\left(\mathcal{L}^{\dagger}\right)^{-1}\left(\lambda^{\dagger,n}\phi^{\dagger,n}\right).
     \end{equation}
     This operation may consist in nested iteration loops: solution of the linear systems and corrections for the non-orthogonal fluxes, non-linearities and explicit terms. The adjoint eigenvalue $\lambda^{\dagger}$  should tend to $\lambda$ for $n \to \infty$.
     \item[2.] Update the eigenvalues using the Rayleigh quotient:
     \begin{equation}
        \label{eq::eigenupdate}
         \lambda^{n+1} = \lambda^n \frac{\left\langle\phi^n \phi^{n+1} \right\rangle}{\left\langle\phi^{n+1} \phi^{n+1} \right\rangle}, \quad \quad
         \lambda^{\dagger,n+1} = \lambda^{\dagger,n} \frac{\left\langle\phi^{\dagger,n} \phi^{\dagger,{n+1}} \right\rangle}{\left\langle\phi^{\dagger,n+1} \phi^{\dagger,{n+1}} \right\rangle}.
     \end{equation}

     While this method is computationally efficient, it results in a convergence rate that is stongly dependent on the initial guess of the eigenvalue. Therefore, we implemented an 'implicit Rayleigh' update which results in an improved convergence rate at higher computational cost:
     \begin{equation}
        \label{eq::eigenupdate2}
         \lambda^{n+1} = \frac{\left\langle\mathcal{L}\phi^{n+1} \phi^{n+1} \right\rangle}{\left\langle\phi^{n+1} \phi^{n+1} \right\rangle}, \quad \quad
         \lambda^{\dagger,n+1} = \frac{\left\langle\mathcal{L}\phi^{\dagger,n+1} \phi^{\dagger,{n+1}} \right\rangle}{\left\langle\phi^{\dagger,n+1} \phi^{\dagger,{n+1}} \right\rangle}.
     \end{equation}
     Notice that the operator $\mathcal{L}$ is computed explicitly after eventual additional terms resulting from non--orthogonal grids have been iterated to convergence.
     
     \item[3.] Normalise $\phi^{n+1}$ and $\phi^{\dagger,n+1}$:
     \begin{equation}
         \phi^{n+1} = \frac{\phi^{n+1} }{\left\langle \phi^{n+1}\right\rangle}, \quad \quad
         \phi^{\dagger,n+1} = \frac{ \phi^{\dagger,n+1} }{\left\langle \phi^{\dagger,n+1} \right\rangle},
     \end{equation}
     Notice that this normalisation is arbitrary and we will later re-normalise $\phi^{\dagger,n+1}$ to be consistent with Eq.~\ref{eq::norm_beta}.
     \item[4.] Check convergence against a number of norms with user-defined tolerances. We choose to test both the residuals for $\phi$ and $\phi^{\dagger}$ defined as:
     \begin{equation}
         \label{eq::residuals_spec}
         \text{res}\left( \phi^{n+1}\right) = \text{max}\left(\frac{|\phi^{n+1}| - |\phi^{n}|}{|\phi^{n}|}\right), \quad \quad \text{res}\left( \phi^{\dagger,n+1}\right) = \text{max}\left(\frac{|\phi^{\dagger,n+1}| - |\phi^{\dagger,n}|}{|\phi^{\dagger,n}|}\right),
     \end{equation}
     Where $\text{max}$ is the maximum and the operator $|\cdot|$ denotes the absolute value. 
     Clearly, the error on the eigenvalues is also a critical metrics to assess convergence:
     \begin{equation}
         (\lambda\text{-error})^{n+1} = \frac{|\lambda^{\dagger,n+1}-\lambda^{n+1}|}{\lambda^{n+1}}.
     \end{equation}
     When all the metrics pass the convergence test ( generally their value should be smaller than $10^{-5}$), the spectral solver exits the loop.
 \end{itemize}
 
 After convergence, the eigenfunctions need to be re-normalised to satisfy $\left\langle \phi \right\rangle = 1$ and $\left\langle \beta \right\rangle = 1$ to be consistent with our formulations. While no action needs to be taken for $\phi$, $\phi^{\dagger}$ is finally re-scaled simply dividing it by $\left\langle \phi^{\dagger}\phi \right\rangle$.

\subsection{Numerical solution of the corrector problem}

Finally, the corrector problem is solved iteratively (with the two nested loops described above). Since $\boldsymbol{\chi}$ is gauge-invariant (i.e., defined up to a constant), from Eq.~\ref{eq::cave_1}, we can make the approximation of $\left\langle c \right\rangle$ of order $\mathcal{O}\left(\eps^2\right)$ by imposing $\left\langle\phi\boldsymbol{\chi}\right\rangle=0$. Thus, at each iteration $n$ we impose:
\begin{equation}
    \boldsymbol{\chi}^{n+1} = \boldsymbol{\chi}^{n+\frac{1}{2}} - \left\langle \phi \boldsymbol{\chi}^{n+\frac{1}{2}} \right\rangle,
\end{equation}
where $\boldsymbol{\chi}^{n+\frac{1}{2}}$ is calculated from the solution of the corrector problem before the gauge scaling.

\subsection{Acceleration of the numerical algorithm}

Power methods can exibith slow convergence, especially when applied to two spectral problems the eigenvalues of which are expected to become equal after a certain number of iterations. Similarly, the rescaling in the corrector problem results in an iterative process. 

In order to accelerate the convergence, we employ \emph{Aitken' acceleration method} (or \emph{Aitken's} $\delta^2$--\emph{process} \citep{numRecipes}). 
This is nothing more than a method for extrapolating the partial sums of a series with approximately geometric convergence, and it is widely used in the solution of stiff fluid--structure interaction problems \citep{Kuttler2008}.

In this work, the eigenfunction, the adjoint function, and the corrector are updated using a relaxation--based implementation of the method of Aitken \citep{Irons1969}:

\begin{equation}
    \phi^{n+1}_{\text{Aitken}} = \phi^{n+1}  + \theta_{\phi}^{n+1}\Delta \phi^{n+1} \,, \quad \Delta \phi^{n} = \phi^n - \phi^{n-1}
    \label{eq::atiken_phi}
\end{equation}
\begin{equation}
    \phi^{\dagger,n+1}_{\text{Aitken}} =\phi^{\dagger,n+1} + \theta_{\phi^{\dagger}}^{n+1}\Delta \phi^{\dagger,n+1} \,, \quad \Delta \phi^{\dagger,n} = \phi^{\dagger,n} - \phi^{\dagger,n-1}
    \label{eq::atiken_phiAdj}
\end{equation}
\begin{equation}
    \boldsymbol{\chi}^{n+1}_{\text{Aitken}} = \boldsymbol{\chi}^{n+1} + \theta_{\boldsymbol{\chi}}^{n+1}\Delta \boldsymbol{\chi}^{n+1} \,, \quad \Delta \boldsymbol{\chi}^{n} = \boldsymbol{\chi}^n - \boldsymbol{\chi}^{n-1}
    \label{eq::atiken_corr}
\end{equation}

where the subscript $\text{Aitken}$ indicates a variable accelerated using Aitken's method. The relaxation factors are computed following Aitken's accelearation method:

\begin{equation}
    \theta_{\phi}^{n+1} = \theta_{\phi}^{n} \frac{\left\langle \Delta \phi^{n+1} (\Delta \phi^{n+1} - \Delta \phi^{n})\right\rangle}{\left\langle |\Delta \phi^{n+1} - \Delta \phi^{n}|\right\rangle}\,, 
    \label{eq::theta_phi}
\end{equation}
\begin{equation}
    \theta_{\phi^{\dagger}}^{n+1} = \theta_{\phi^{\dagger}}^{n} \frac{\left\langle \Delta \phi^{\dagger,n+1} (\Delta \phi^{\dagger,n+1} - \Delta \phi^{\dagger,n})\right\rangle}{\left\langle |\Delta \phi^{\dagger,n+1} - \Delta \phi^{\dagger,n}|\right\rangle}\,,    
    \label{eq::theta_phiAdj}
\end{equation}
\begin{equation}
    \theta_{\boldsymbol{\chi}}^{n+1} = \theta_{\boldsymbol{\chi}}^{n} \frac{\left\langle \Delta \boldsymbol{\chi}^{n+1} \cdot (\Delta \boldsymbol{\chi}^{n+1} - \Delta \boldsymbol{\chi}^{n})\right\rangle}{\left\langle |\Delta \boldsymbol{\chi}^{n+1} - \Delta \boldsymbol{\chi}^{n}|\right\rangle}\,.  
    \label{eq::theta_corr}
\end{equation}

Notice that in our algorithm the acceleration is performed uniformly on all the cells (i.e., the relaxation factors are just numbers and not fields).

\section{Numerical results}
\label{S:res}

\subsection{Verification}

\begin{figure}[htbp]
    \centering
    \includegraphics[width=\linewidth]{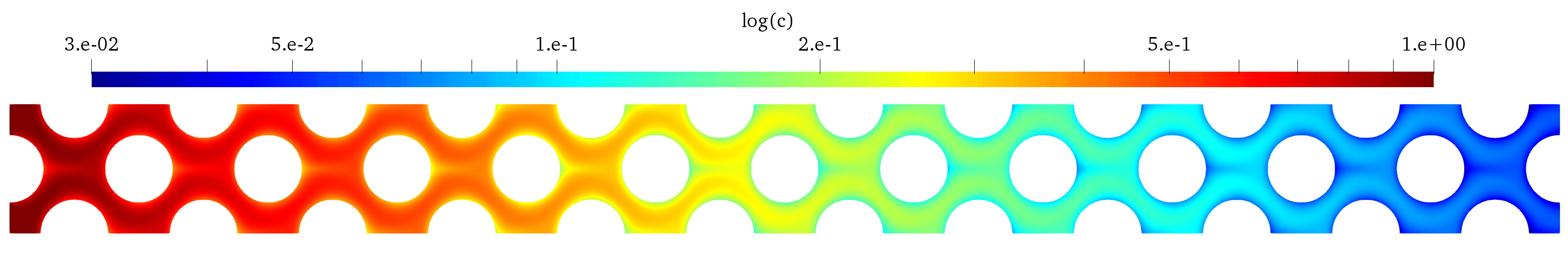}
    \caption{Example results from a two-dimensional pore-scale simulation. For clarity, we show only a fraction of the FCC cells composing the computational domain.}
    \label{fig::poreScale}
\end{figure}

We verify both our code and the upscaling methodology by direct comparison with spatial averaged data from fully resolved pore-scale simulations. Flow and scalar transport are solved in two dimensions for an array of 26 face-centred-cubic (FCC) cells (see Fig.~\ref{fig::poreScale}) using the OpenFOAM\textsuperscript{\textregistered} native solvers \emph{simpleFoam} (classic Navier-Stokes solver employing the SIMPLE algorithm for pressure-velocity coupling ) and \emph{scalarTransportFoam} (standard advection-diffusion equation corresponding to Eq.~\ref{eq::ADEdimensional}). To implement the Robin boundary condition we follow \cite{boccardo2018robust}.

When solving Eq.~\ref{eq::ADEdimensional} at the pore-scale, we provide the following external boundary conditions:
\begin{equation}
    \label{eq::porescale_BCs}
    c(x = 0, y) = 1, \quad \left.\frac{\partial c}{\partial x} \right|_{x=L} = 0,
\end{equation}
where $x$ is the axial direction of the cell array and $L$ is the domain length.
We ensure the flow is in Stokes (viscous) regime by imposing a value of the Reynolds number $\RE < 10^{-3}$ everywhere. 
Steady-state results are then averaged over each cell and compared against predictions obtained from the ordinary differential equation (ODE):
\begin{equation}
    \label{eq::ODE_chebfun}
    \frac{d}{dx}\left( V^{\star}_{x}  \langle c \rangle - \mathcal{D}^{\text{eff}}_{xx}\frac{d  \langle c \rangle }{dx}\right) = -\lambda  \langle c \rangle \, ,
\end{equation}

 where $V^{\star}_x$ is the effective velocity in $x$ and $\mathcal{D}^{\text{eff}}_{xx}$ is the axial component of the effective diffusivity tensor.
Equation~\ref{eq::ODE_chebfun} is solved to spectral accuracy using the MATLAB\textsuperscript{\textregistered} package \emph{Chebfun} \cite{chebfun}. 
 
 Choosing appropriate boundary conditions for Eq.~\ref{eq::ODE_chebfun} is not trivial, since we do not know the value of $\langle c \rangle$ at $x=0$. However, since our only objective is to evaluate the accuracy of this method, we  just impose:
 \begin{equation}
    \label{eq::ODE_BCs}
    \langle c \rangle (x = 0) = 1, \quad \left.\frac{d \langle c \rangle}{d x} \right|_{x=L} = 0,
\end{equation}
and compare the results against fully developed (i.e., far from the inlet) pore-scale simulations with an appropriate rescaling. Therefore, all comparisons will be made dividing all values of $\langle c \rangle$ obtained from pore-scale simulations with the value of $\langle c \rangle$ at the 10$^{th}$ FCC cell, where the profile of $c$ is well developed for all the simulations.

\begin{figure}[htbp]
    \centering
    \begin{subfigure}[htbp]{0.5\textwidth}
        \centering
        \includegraphics[width=0.9\textwidth]{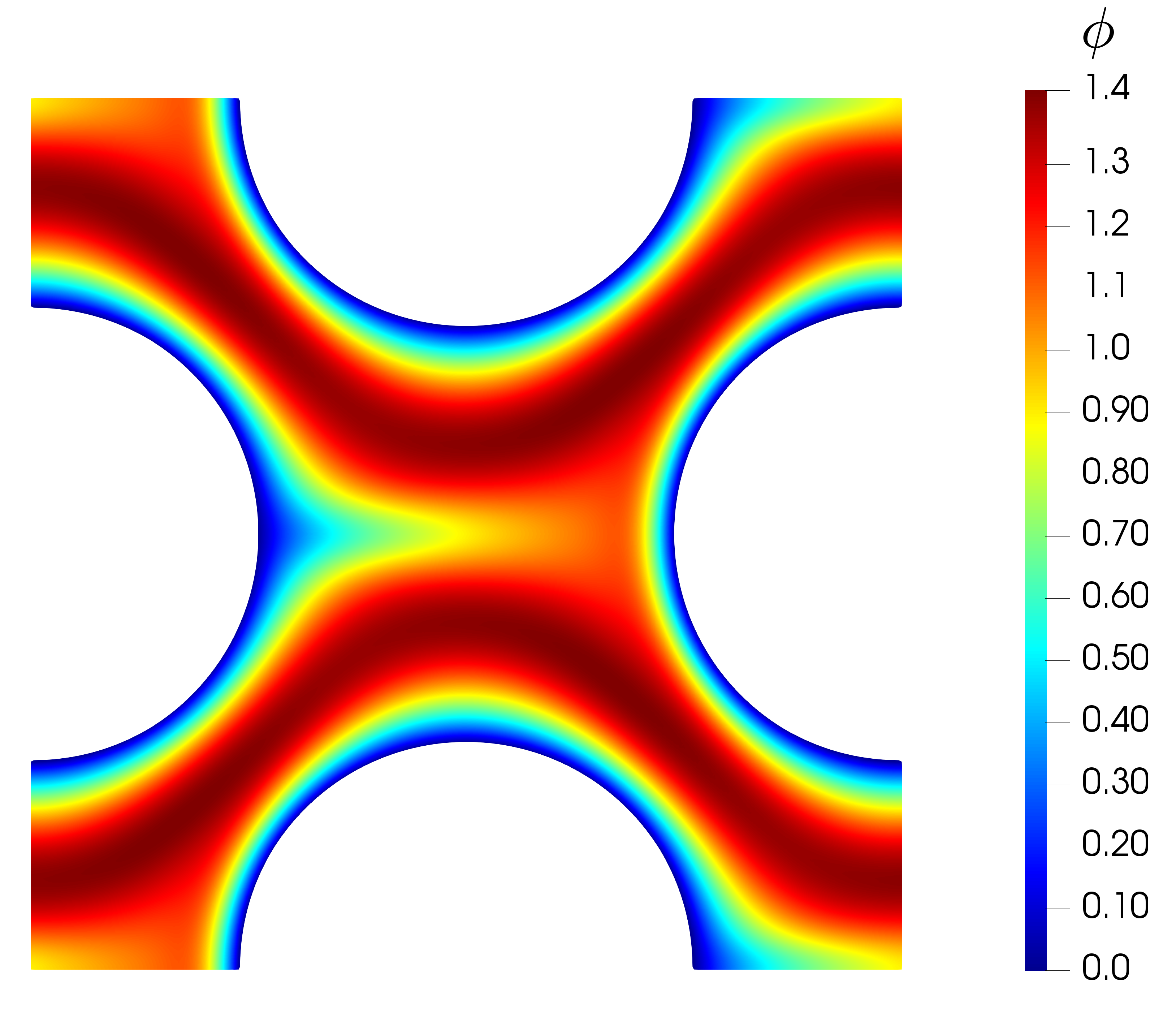}
        \caption{}
        \label{subfig::phi}
    \end{subfigure}%
    ~ 
    \begin{subfigure}[htbp]{0.5\textwidth}
        \centering
        \includegraphics[width=0.9\textwidth]{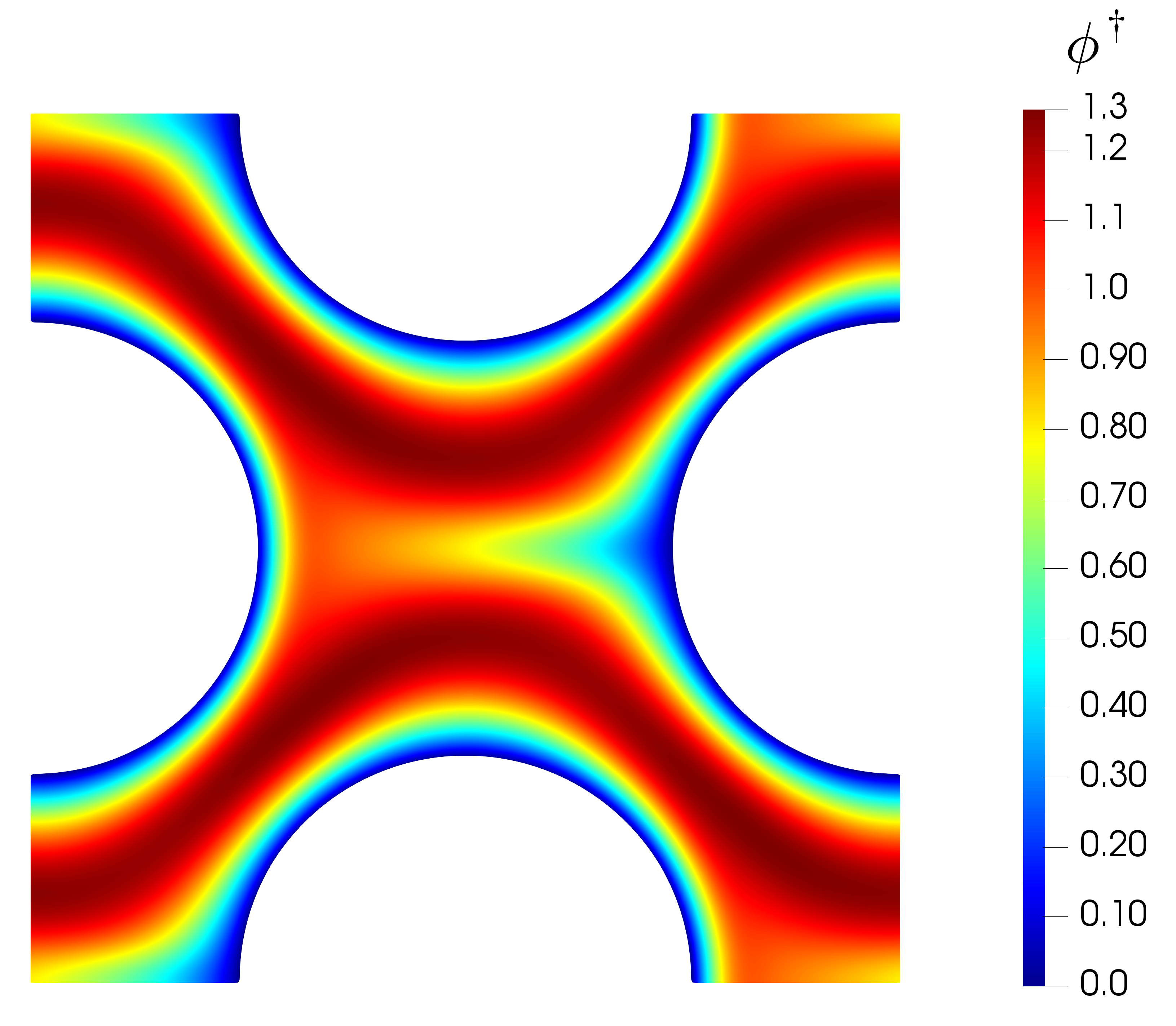}
        \caption{}
        \label{subfig::phiAdj}
    \end{subfigure}%
    \\
    \begin{subfigure}[htbp]{0.5\textwidth}
        \centering
        \includegraphics[width=0.9\textwidth]{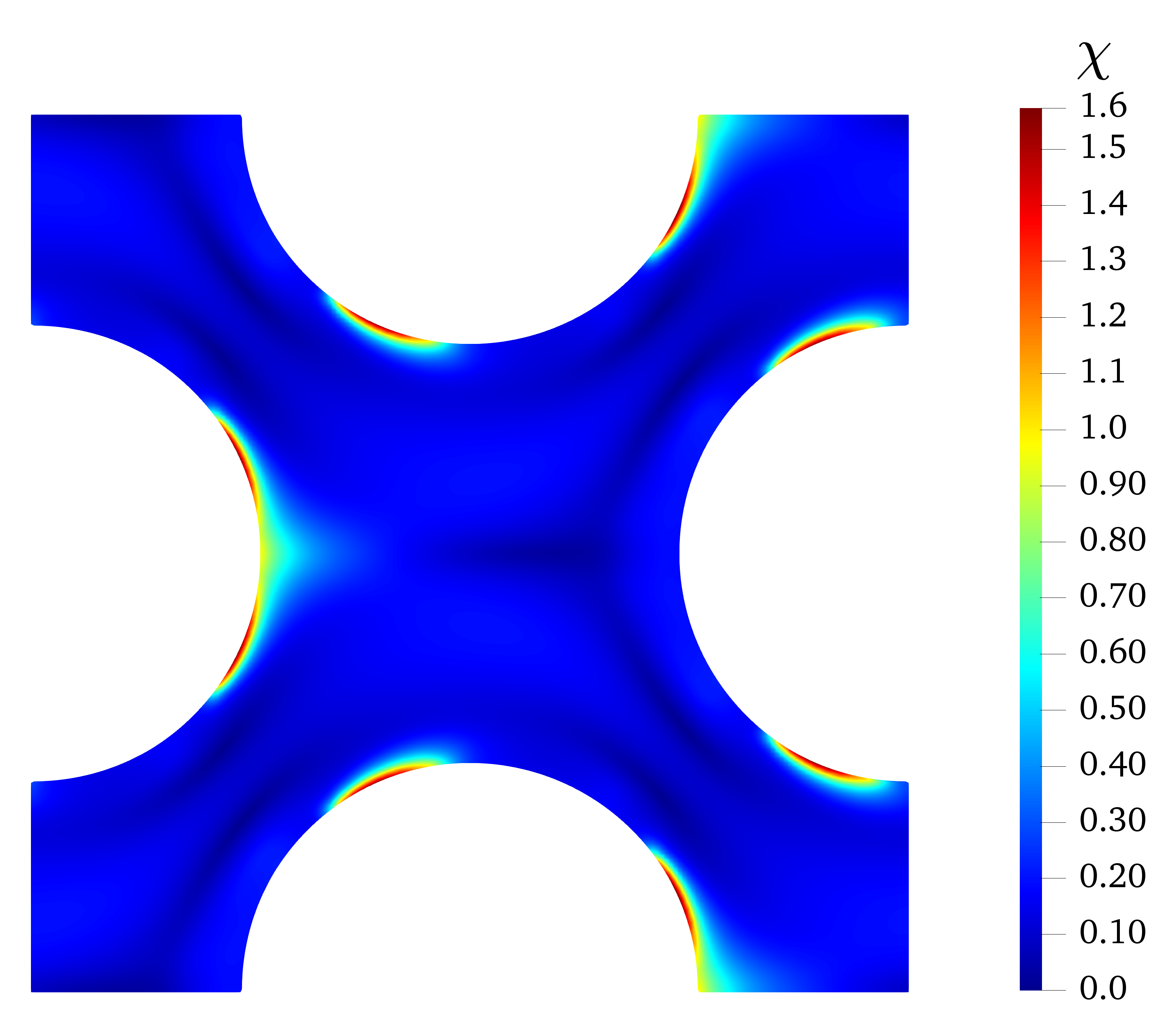}
        \caption{}
        \label{subfig::chi}
    \end{subfigure}%
    \begin{subfigure}[htbp]{0.5\textwidth}
        \centering
        \includegraphics[width=0.9\textwidth]{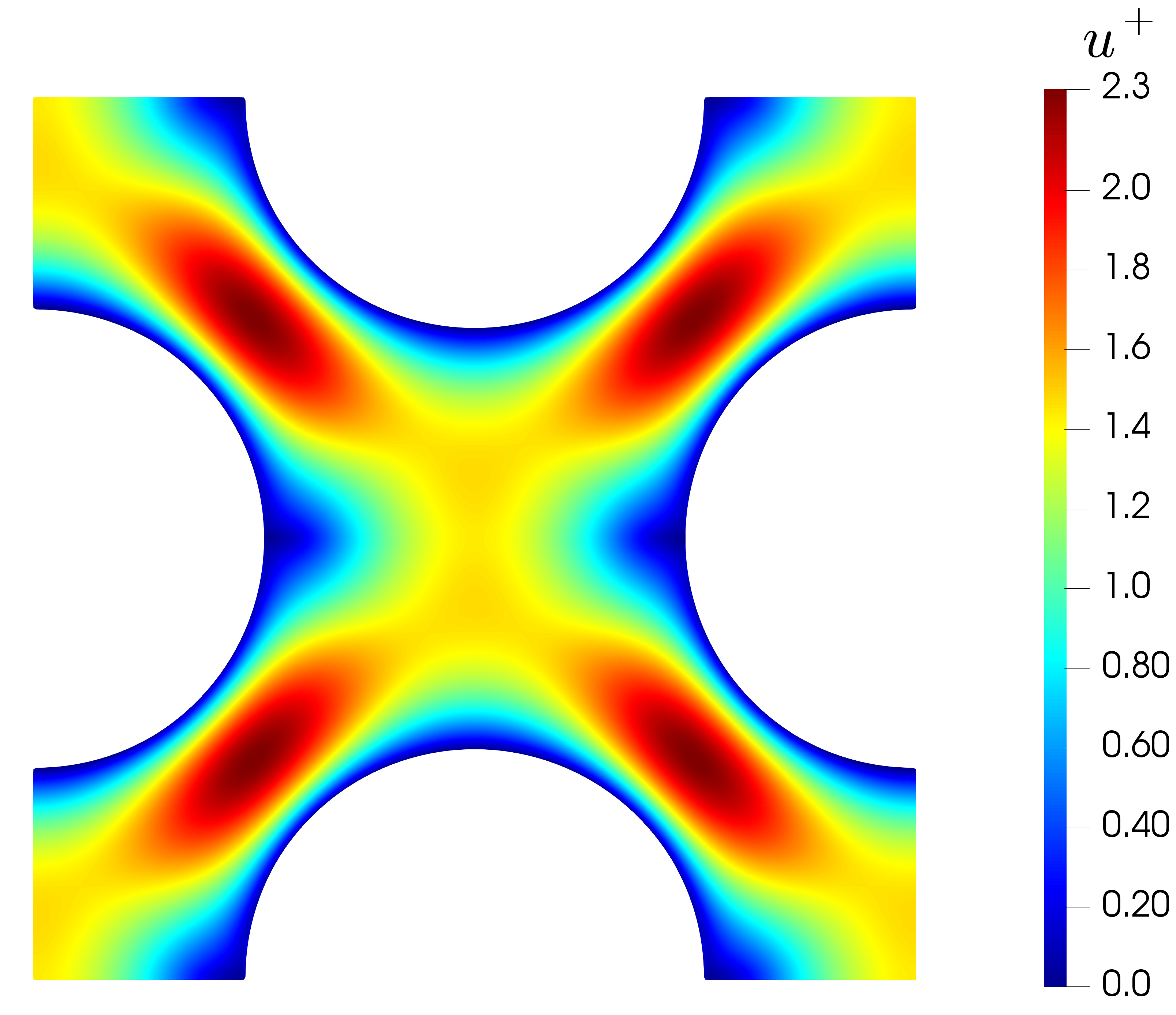}
        \caption{}
        \label{subfig::uplus}
    \end{subfigure}%
    \caption{Results from the spectral and first order corrector solver for the FCC cell in the case $\eps^{-1} \PE = 100$, $\epsilon = 0.7$ and $\eps^{-1} \DAII = 962$. $u^{+}$ is the magnitude of $\bm{u}^{+}$ and $\chi$ is the magnitude of $\boldsymbol{\chi}$.}
    \label{fig::cellRes}
\end{figure}
We compute the coefficients in Eq.~\ref{eq::ODE_chebfun} using our novel solver from a single FCC cell. Fig.~\ref{fig::cellRes} shows the relevant fields arising from the solution of the cell problem. Notice that the eigenfunction $\phi$ in Fig.~\ref{subfig::phi} and the adjoint $\phi^{\dagger}$ in Fig.~\ref{subfig::phiAdj} only differ for the direction of the advective component and for the scaling factor, as expected from inspection of their governing equation and the choice made for their normalisation.

 \begin{figure}[htbp]

    \centering
    \begin{subfigure}[htbp]{0.5\textwidth}
        \centering
        \includegraphics[width=0.9\textwidth]{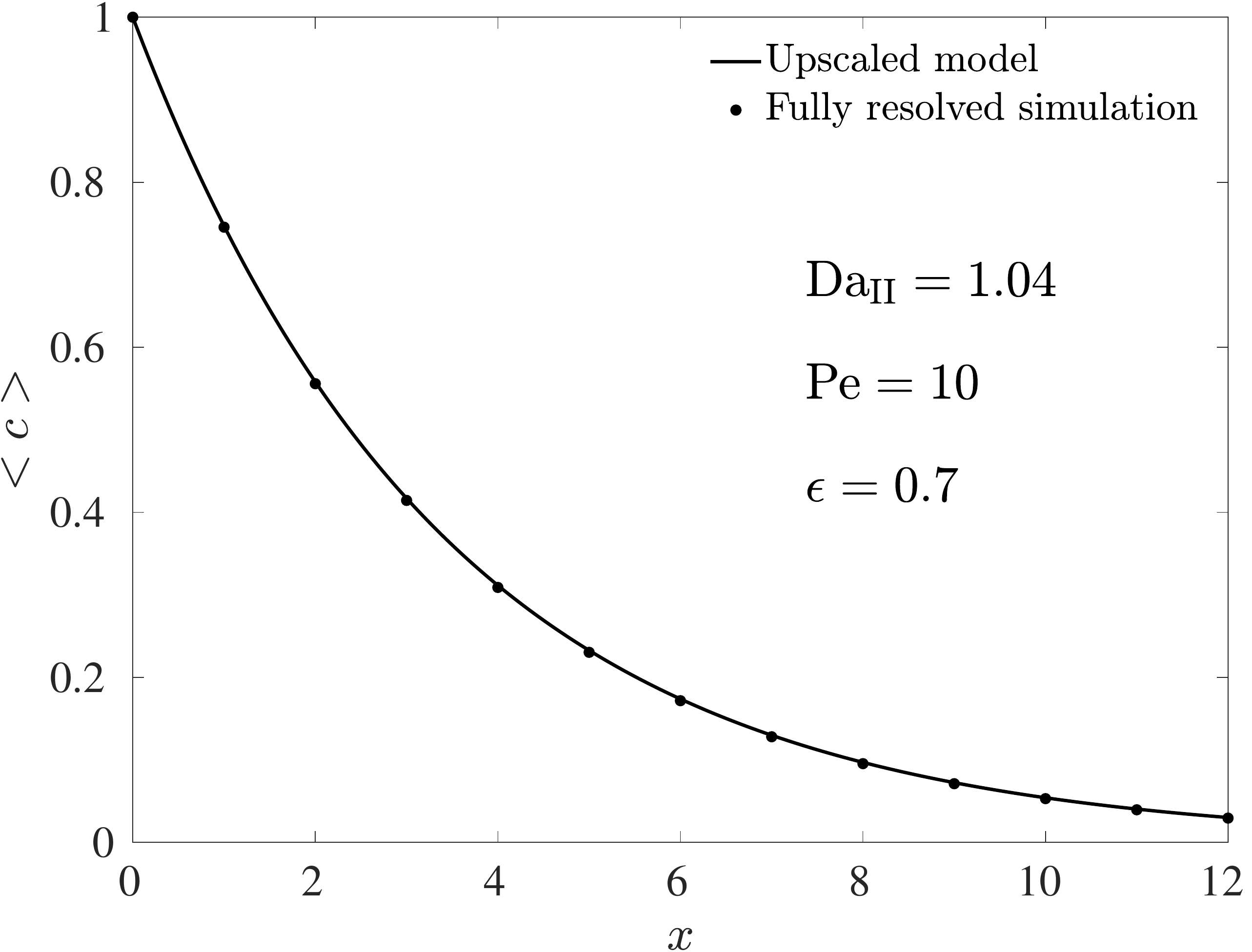}
        \caption{}
        \label{subfig::val1}
    \end{subfigure}%
    ~ 
    \begin{subfigure}[htbp]{0.5\textwidth}
        \centering
        \includegraphics[width=0.9\textwidth]{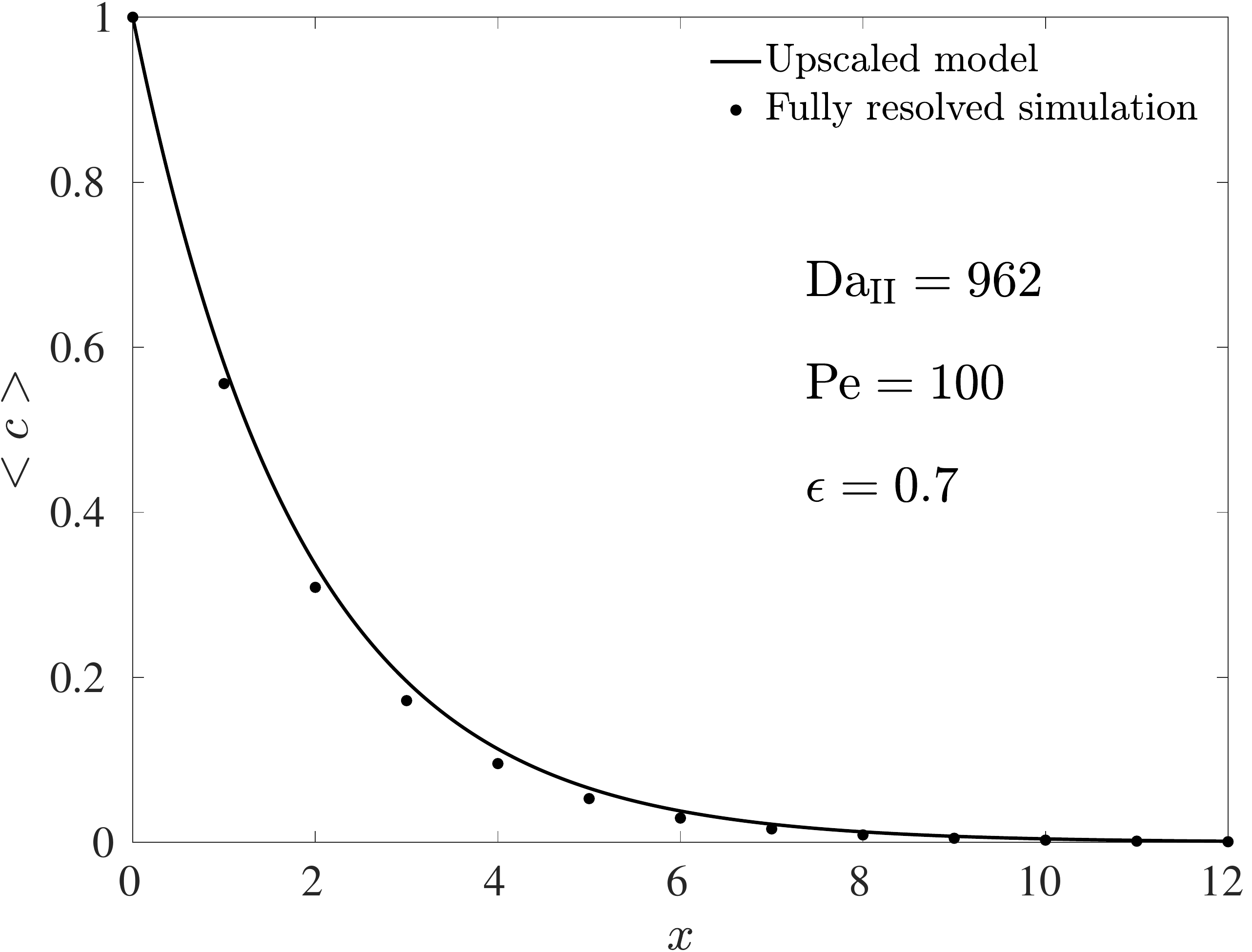}
        \caption{}
        \label{subfig::val2}
    \end{subfigure}%
    \\
    \begin{subfigure}[htbp]{0.5\textwidth}
        \centering
        \includegraphics[width=0.9\textwidth]{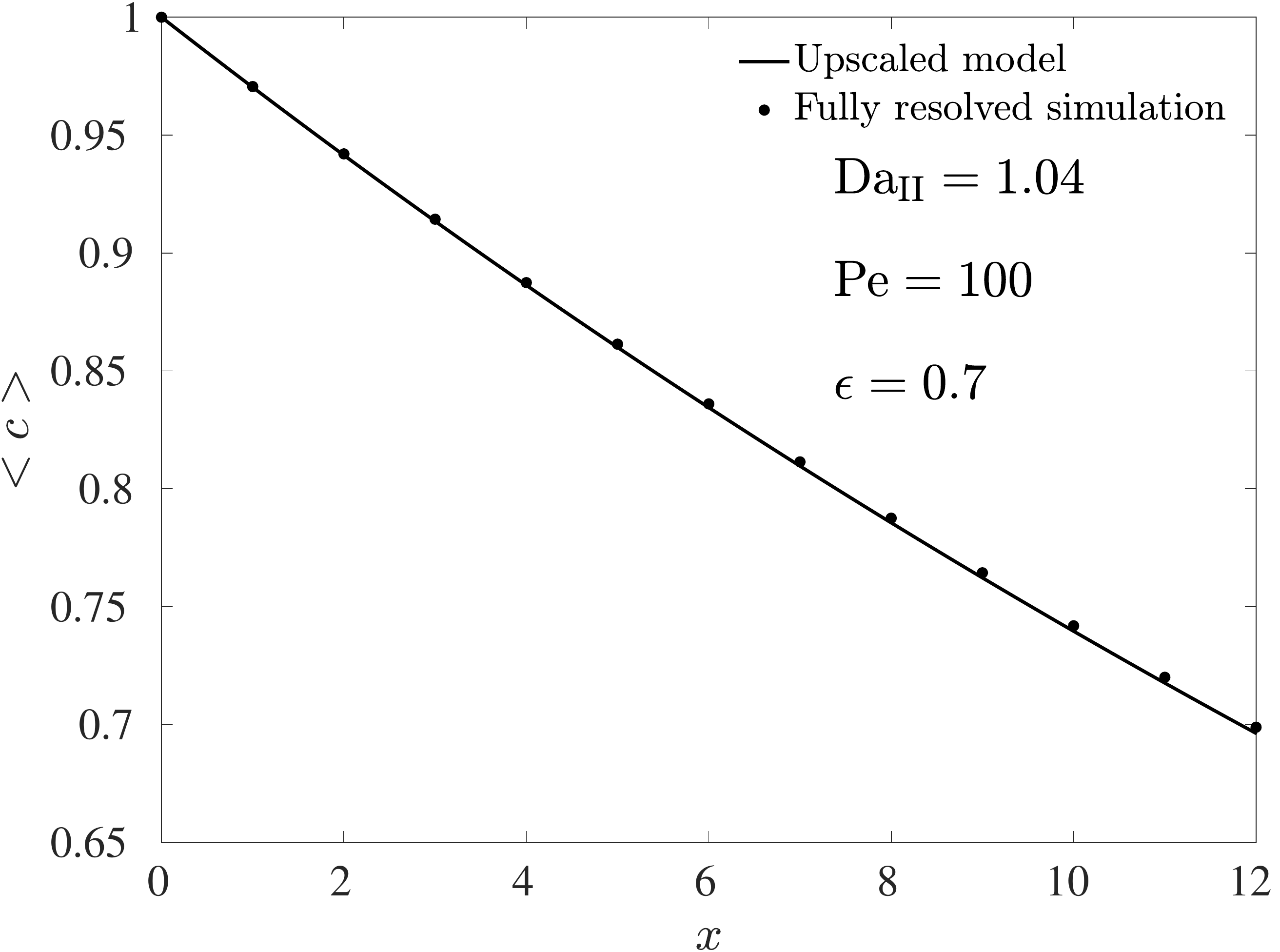}
        \caption{}
        \label{subfig::val3}
    \end{subfigure}%
    ~ 
    \begin{subfigure}[htbp]{0.5\textwidth}
        \centering
        \includegraphics[width=0.9\textwidth]{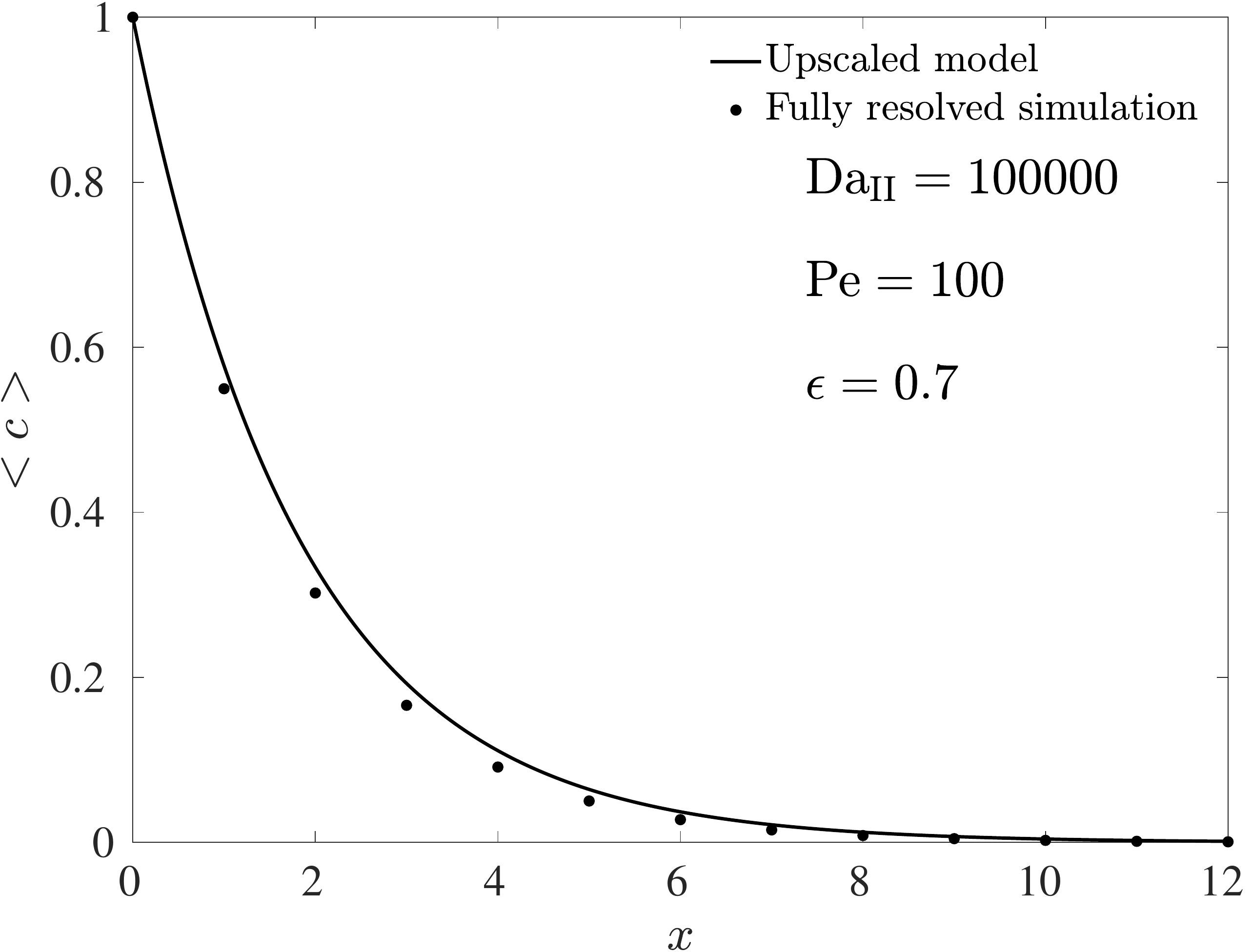}
        \caption{}
        \label{subfig::val4}
    \end{subfigure}%
    \caption{Comparison of results from Eq.~\ref{eq::ODE_chebfun} against resolved pore-scale simulations. Here, $x$ is scaled with the lenght of an FCC cell such that the value of $x$ corresponds to the number of FCC cells.}
    \label{fig::validation}
\end{figure}

Results from two pore-scale simulations are compared with solutions of Eq.~\ref{eq::ODE_chebfun} in Fig~\ref{fig::validation}. Overall, the upscaling method provides excellent results, with little deviations from the pore-scale simulations. Notice that such agreement was obtained by rescaling the cell average concentration from the pore-scale simulations by an appropriate reference value (i.e., the value of an FCC cell in the asymptotic regime). Using this approach, we were able to test the accuracy of the method without a complete knowledge of the external boundary conditions.

\subsection{Parametric study - FCC array}

\begin{figure}[htbp]

    \centering
    \begin{subfigure}[htbp]{0.5\textwidth}
        \centering
        \includegraphics[width=0.9\textwidth]{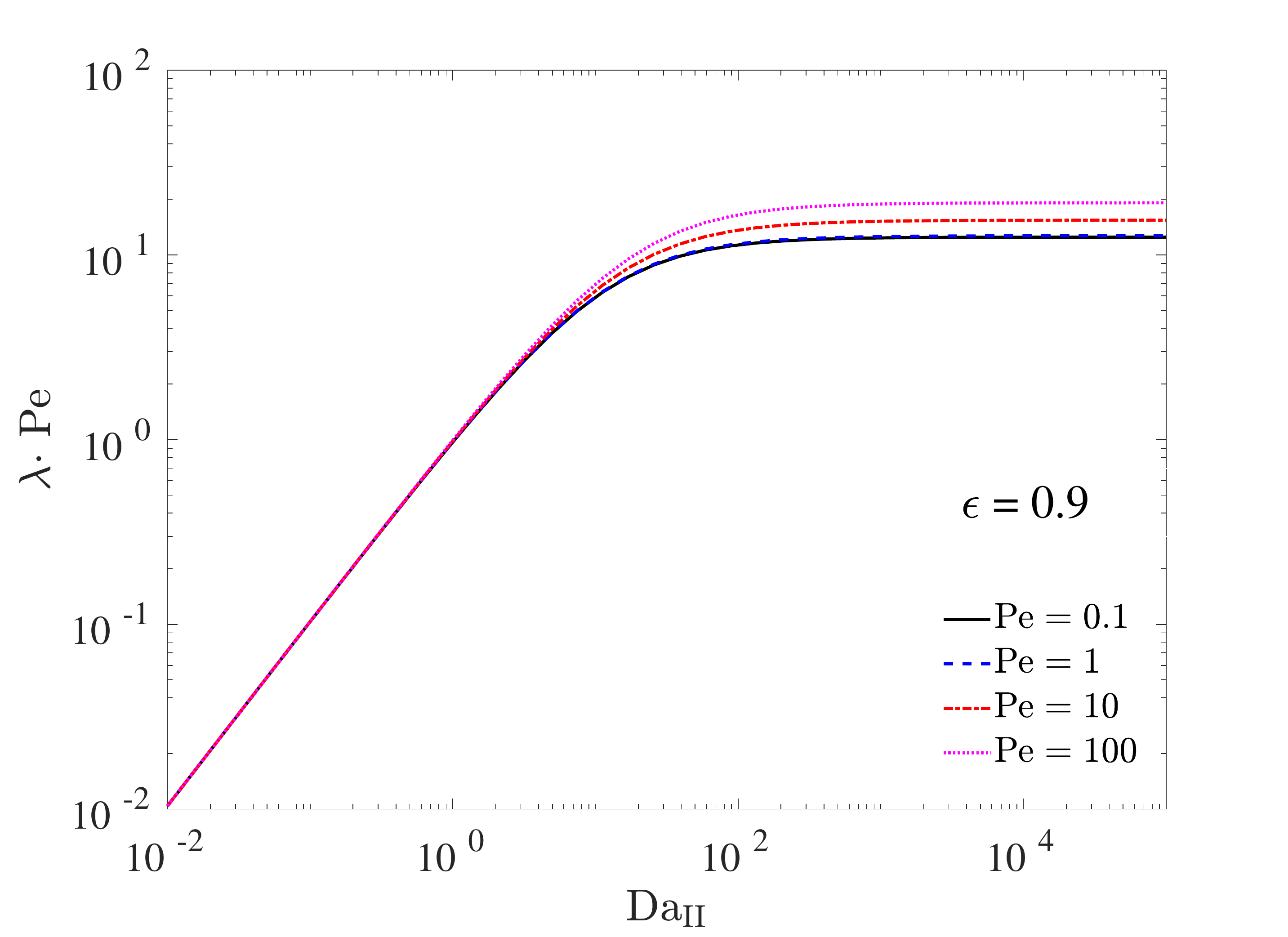}
        \caption{}
        \label{subfig::eps01_lambda}
    \end{subfigure}%
    ~ 
    \begin{subfigure}[htbp]{0.5\textwidth}
        \centering
        \includegraphics[width=0.9\textwidth]{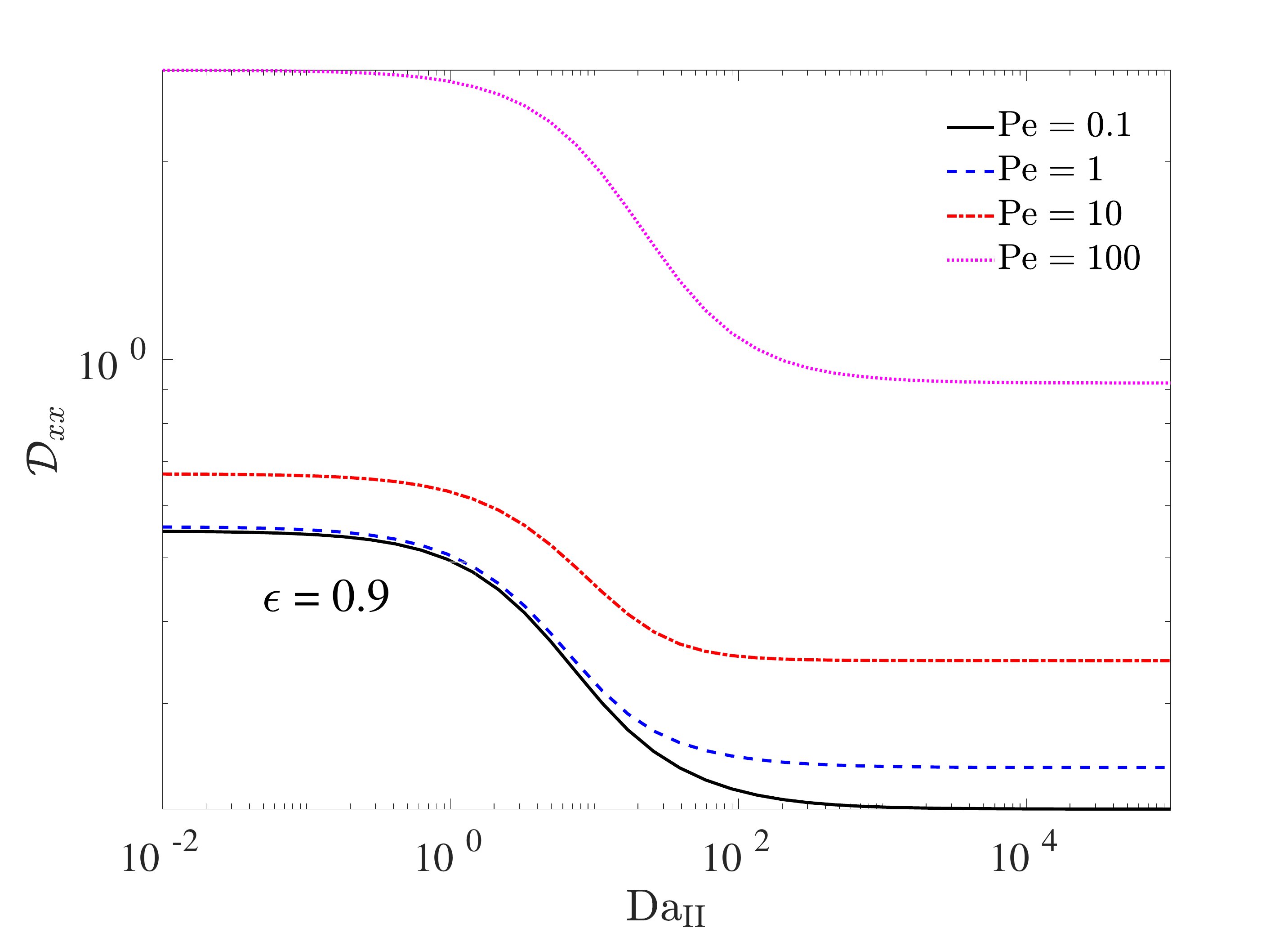}
        \caption{}
        \label{subfig::eps01_Dx}
    \end{subfigure}%
    \\
    \begin{subfigure}[htbp]{0.5\textwidth}
        \centering
        \includegraphics[width=0.9\textwidth]{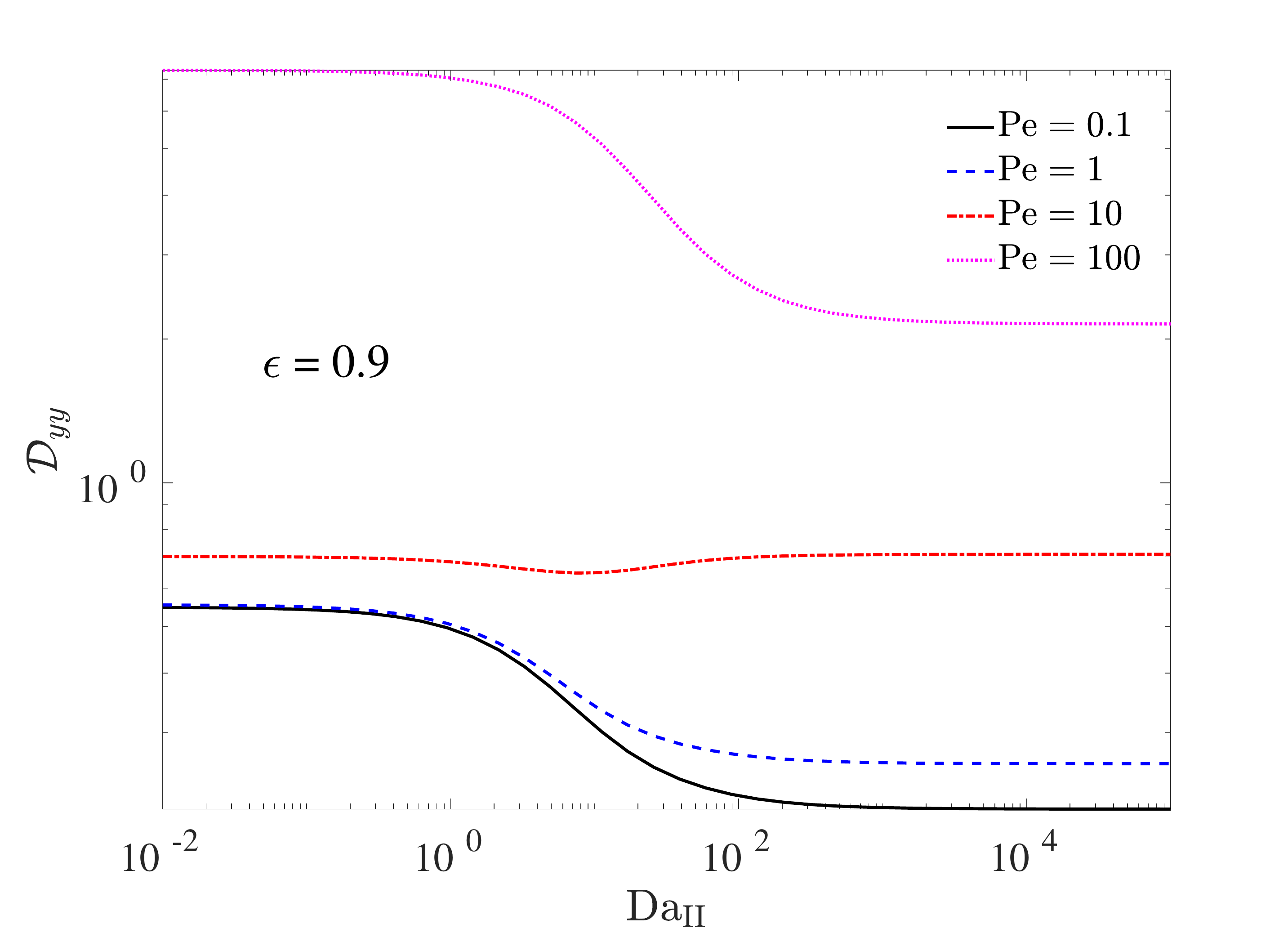}
        \caption{}
        \label{subfig::eps01_Dy}
    \end{subfigure}%
    ~ 
    \begin{subfigure}[htbp]{0.5\textwidth}
        \centering
        \includegraphics[width=0.9\textwidth]{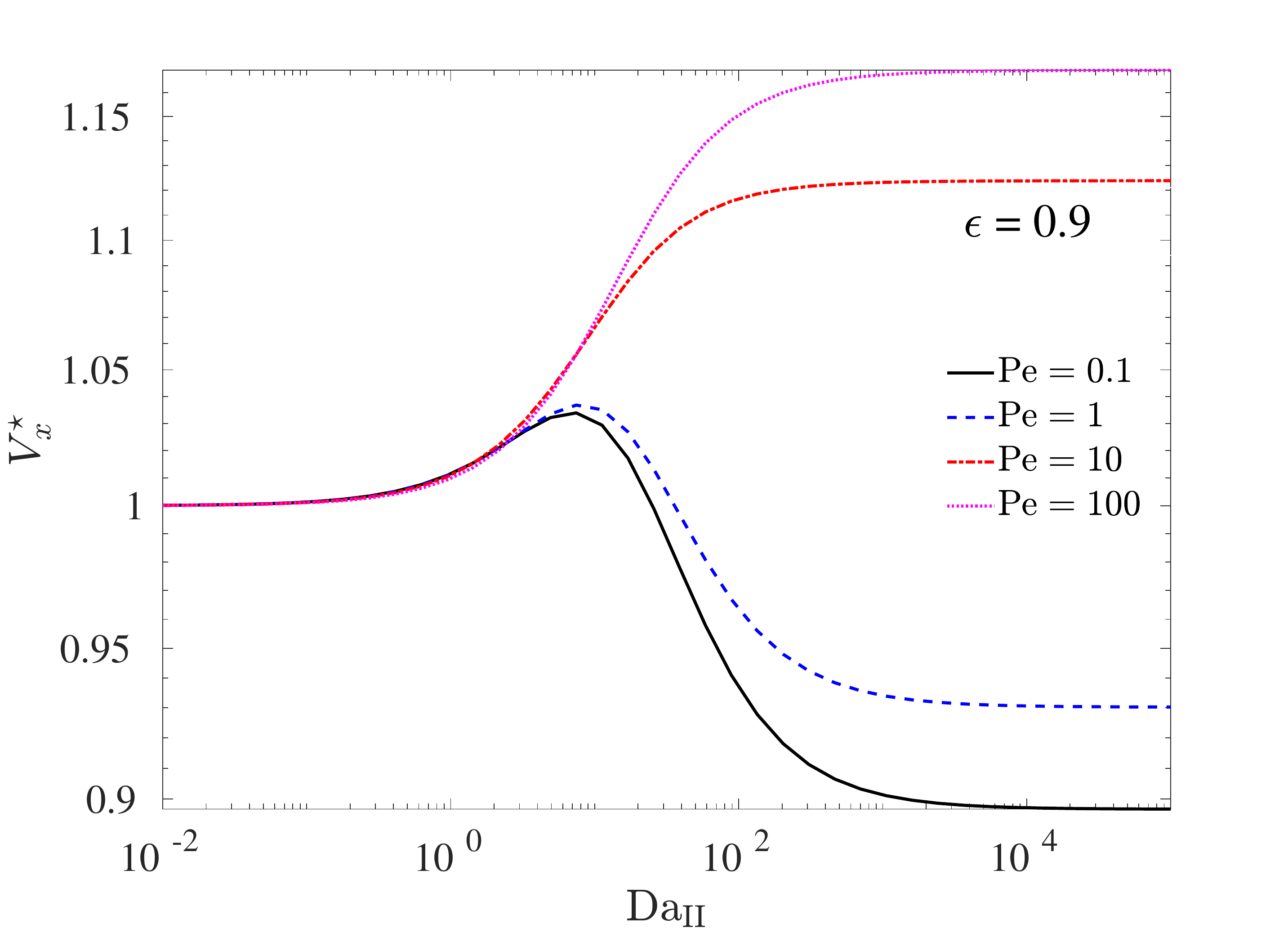}
        \caption{}
        \label{subfig:eps01_U}
    \end{subfigure}%
    \caption{Scaled effective parameters as a function of the microscopic parameters for $\epsilon=0.9$.}
    \label{fig::eps01}
\end{figure}

\begin{figure}[htbp]

    \centering
    \begin{subfigure}[htbp]{0.5\textwidth}
        \centering
        \includegraphics[width=0.9\textwidth]{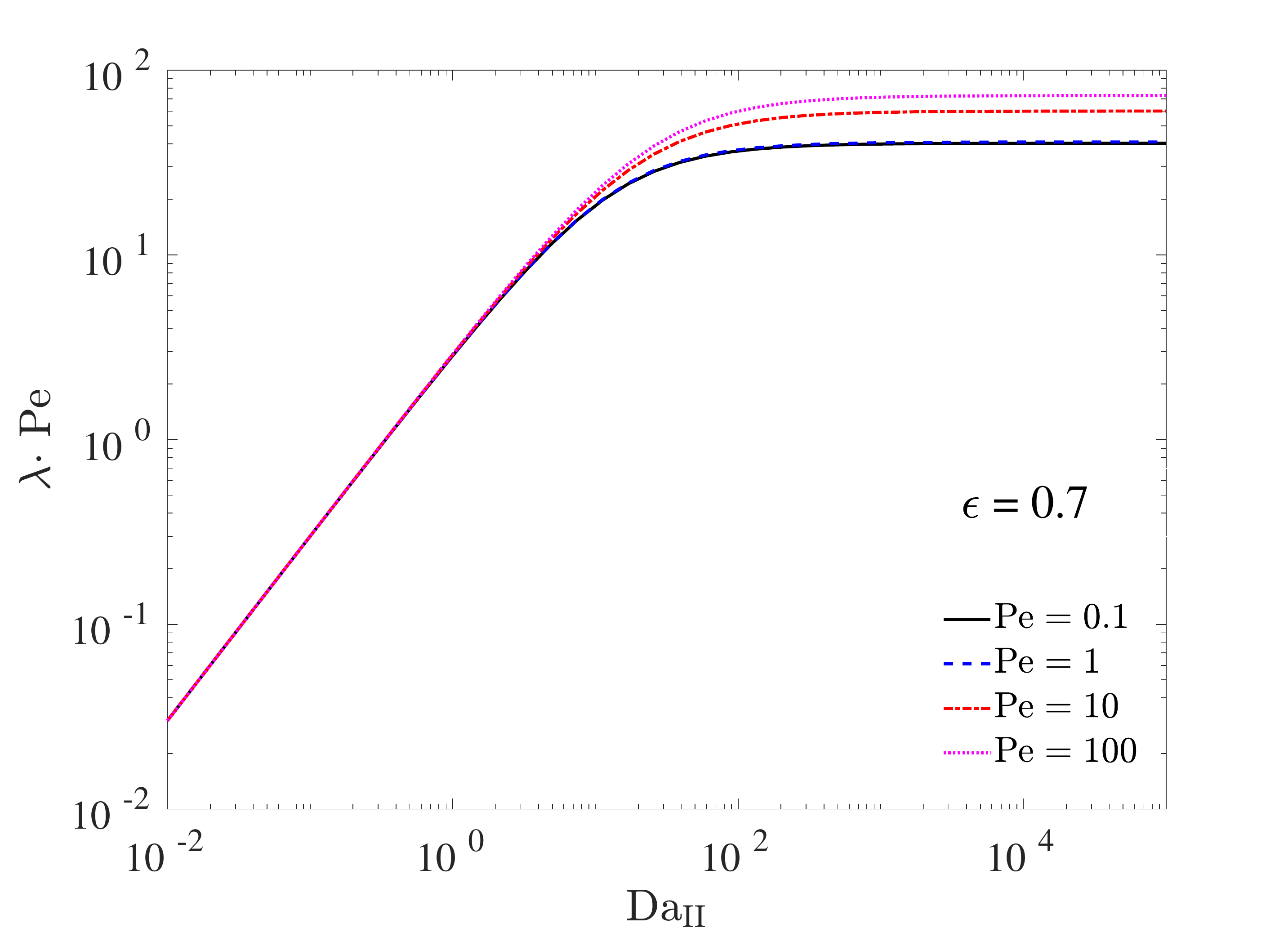}
        \caption{}
        \label{subfig::eps07_lambda}
    \end{subfigure}%
    ~ 
    \begin{subfigure}[htbp]{0.5\textwidth}
        \centering
        \includegraphics[width=0.9\textwidth]{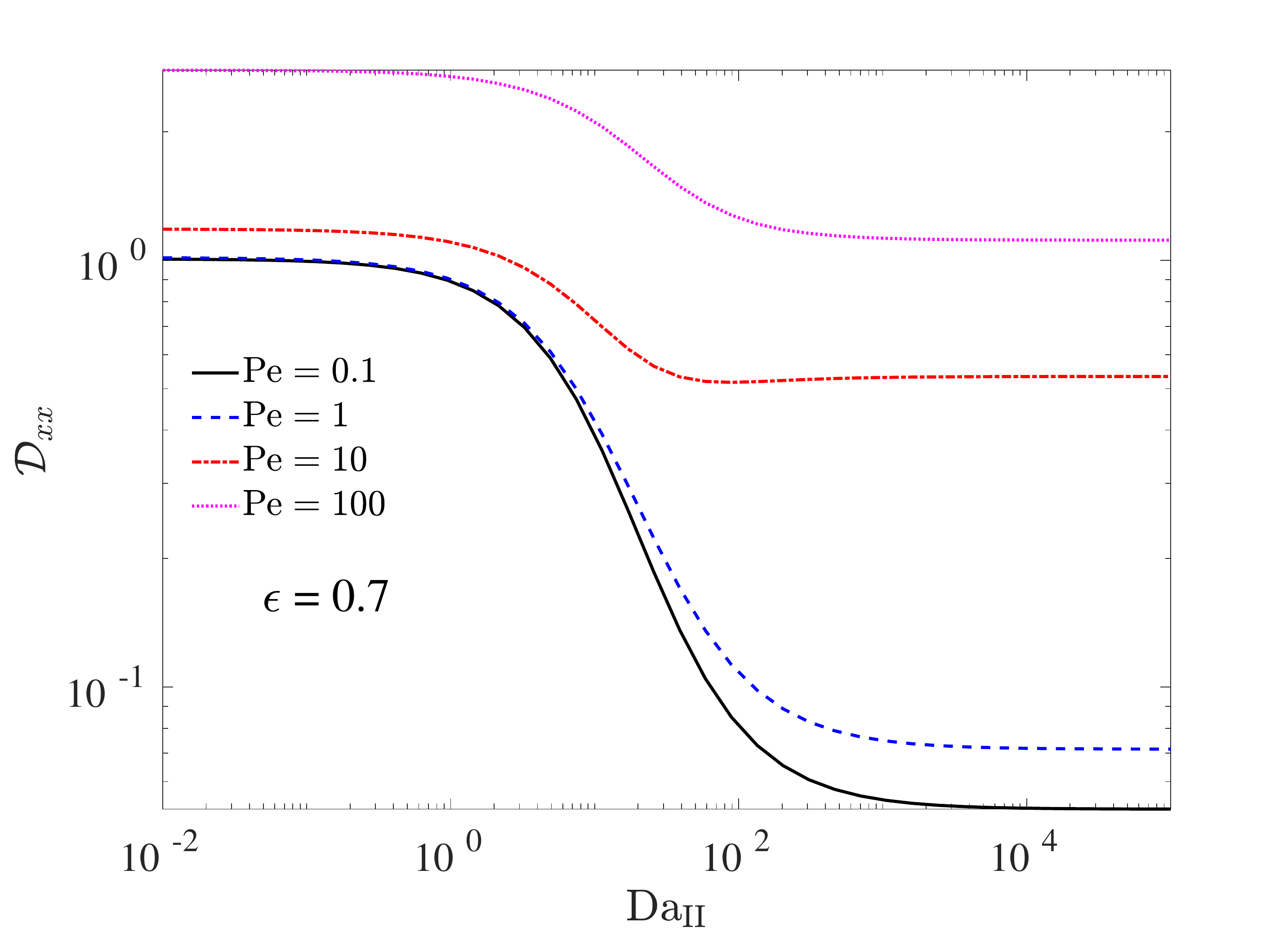}
        \caption{}
        \label{subfig::eps07_Dx}
    \end{subfigure}%
    \\
    \begin{subfigure}[htbp]{0.5\textwidth}
        \centering
        \includegraphics[width=0.9\textwidth]{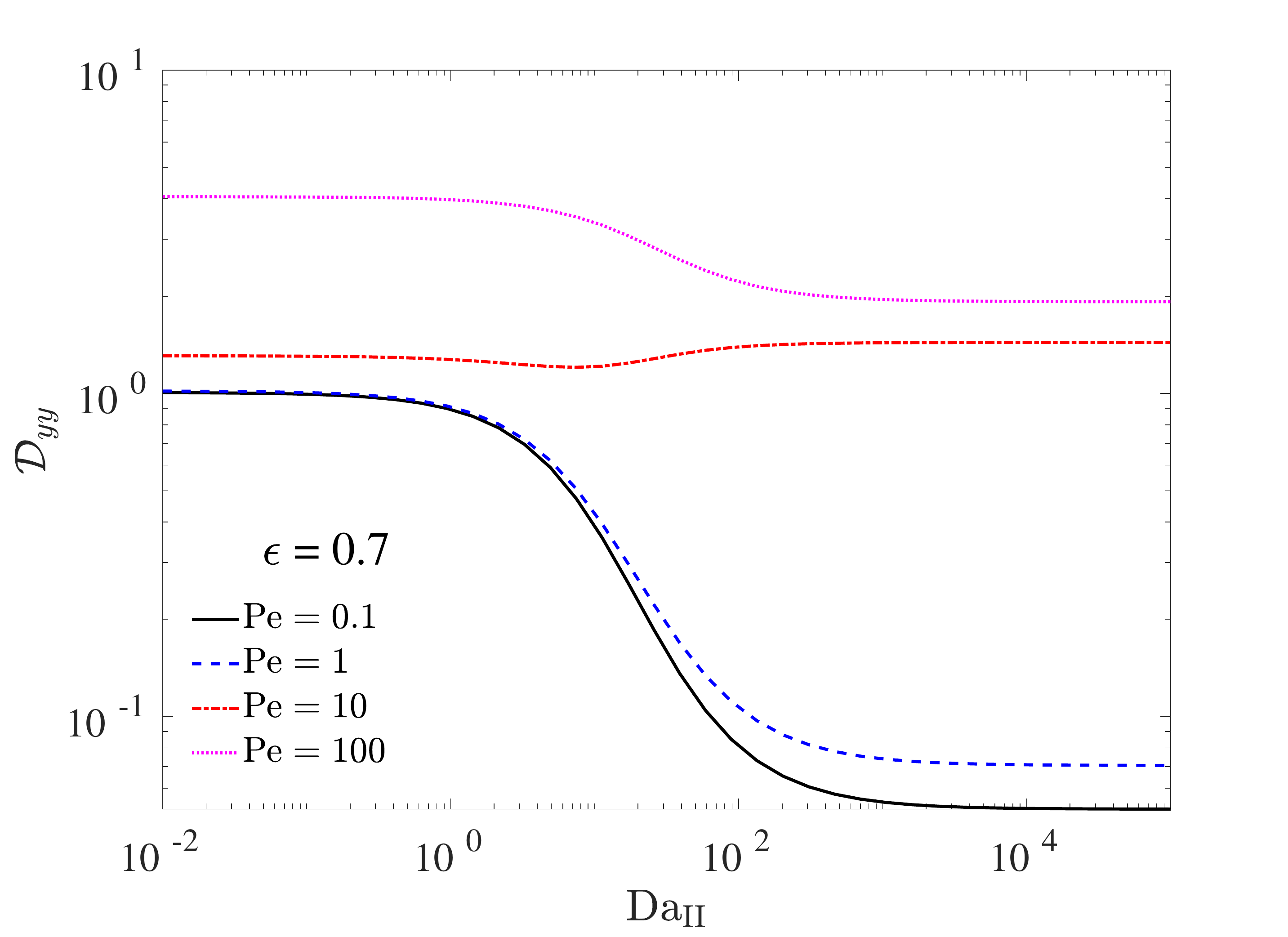}
        \caption{}
        \label{subfig::eps07_Dy}
    \end{subfigure}%
    ~ 
    \begin{subfigure}[htbp]{0.5\textwidth}
        \centering
        \includegraphics[width=0.9\textwidth]{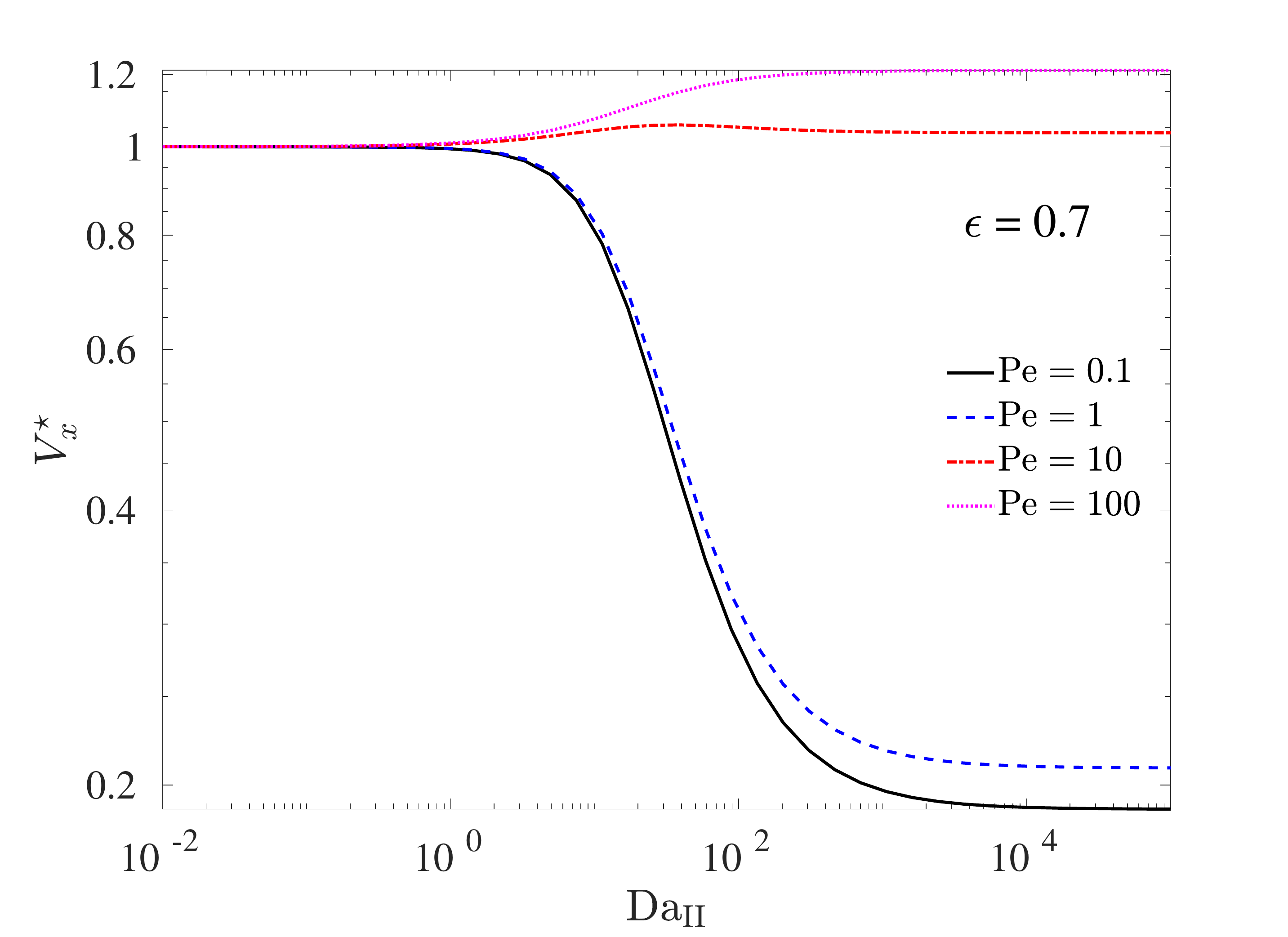}
        \caption{}
        \label{subfig:eps07_U}
    \end{subfigure}%
    \caption{Scaled effective parameters as a function of the microscopic parameters for $\epsilon=0.7$.}
    \label{fig::eps07}
\end{figure}

\begin{figure}[htbp]

    \centering
    \begin{subfigure}[htbp]{0.5\textwidth}
        \centering
        \includegraphics[width=0.9\textwidth]{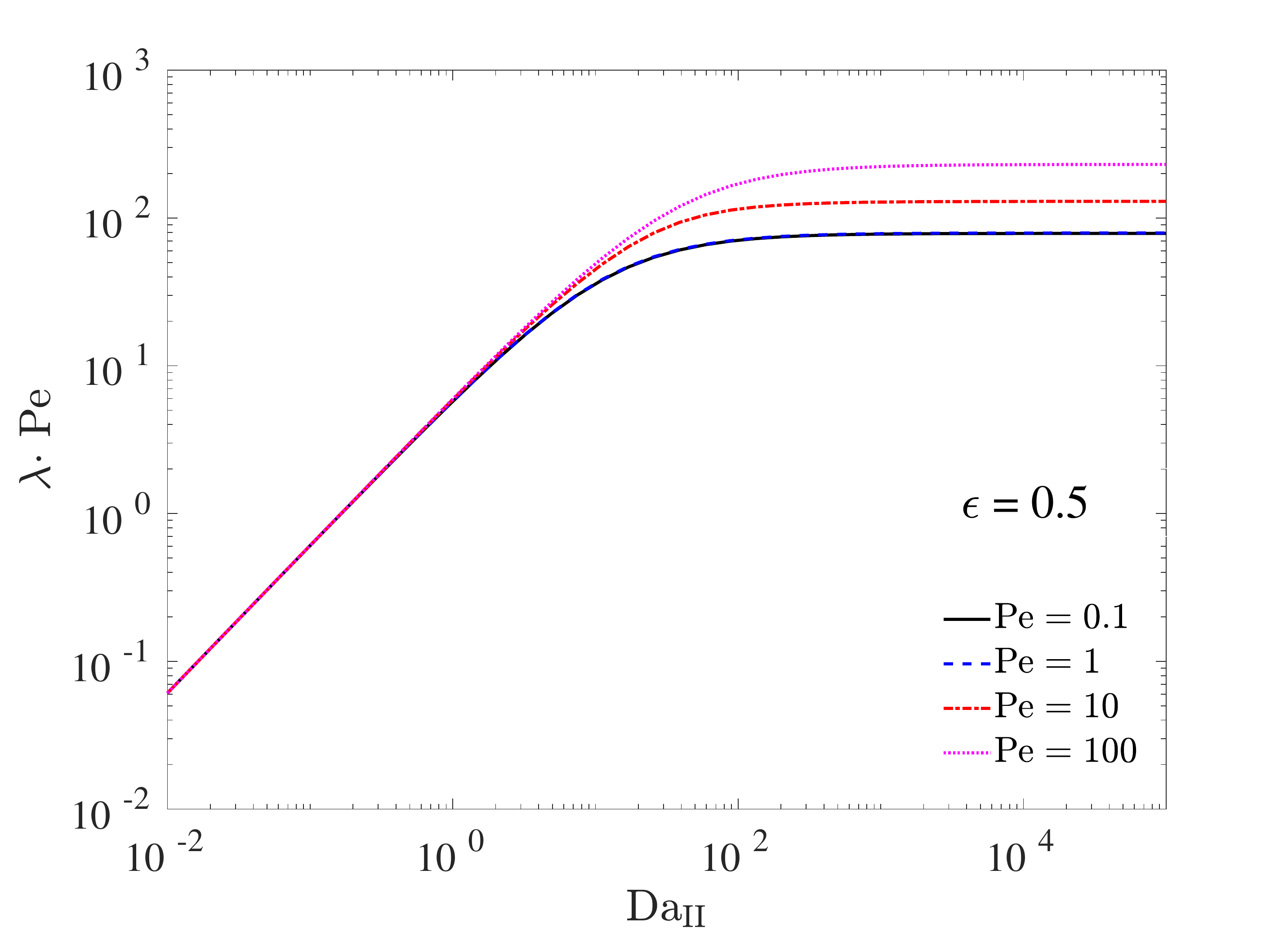}
        \caption{}
        \label{subfig::eps05_lambda}
    \end{subfigure}%
    ~ 
    \begin{subfigure}[htbp]{0.5\textwidth}
        \centering
        \includegraphics[width=0.9\textwidth]{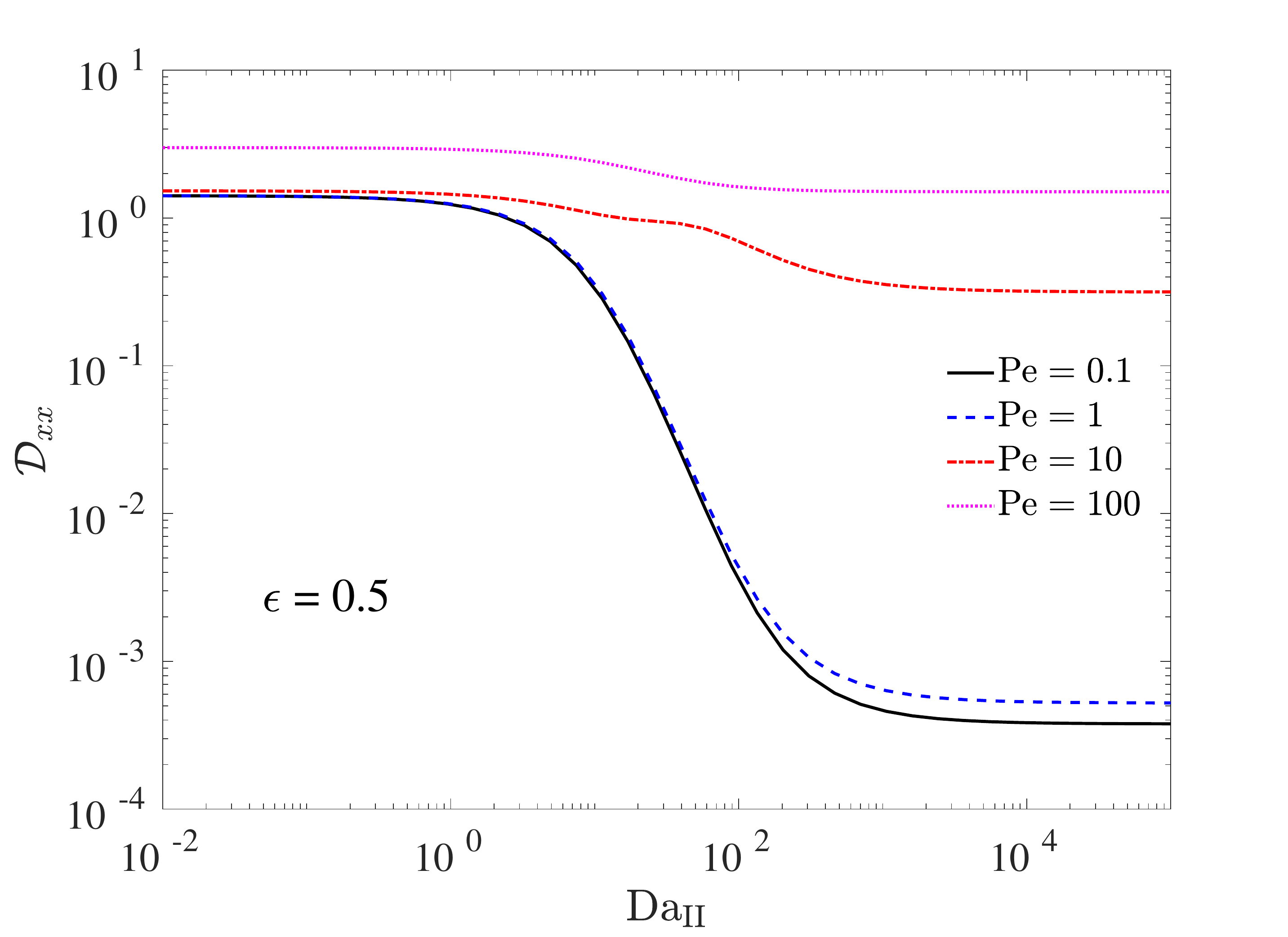}
        \caption{}
        \label{subfig::eps05_Dx}
    \end{subfigure}%
    \\
    \begin{subfigure}[htbp]{0.5\textwidth}
        \centering
        \includegraphics[width=0.9\textwidth]{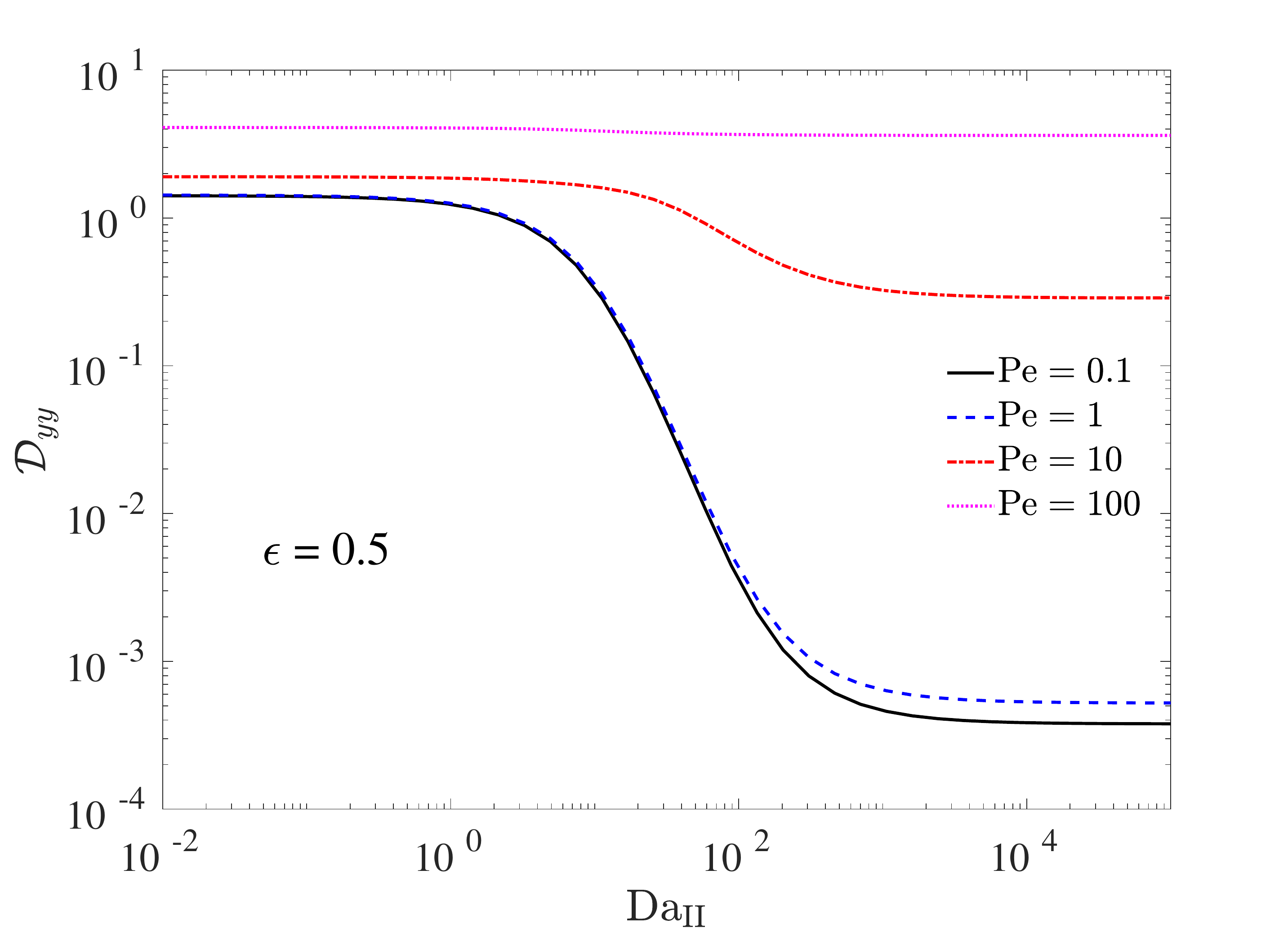}
        \caption{}
        \label{subfig::eps05_Dy}
    \end{subfigure}%
    ~ 
    \begin{subfigure}[htbp]{0.5\textwidth}
        \centering
        \includegraphics[width=0.9\textwidth]{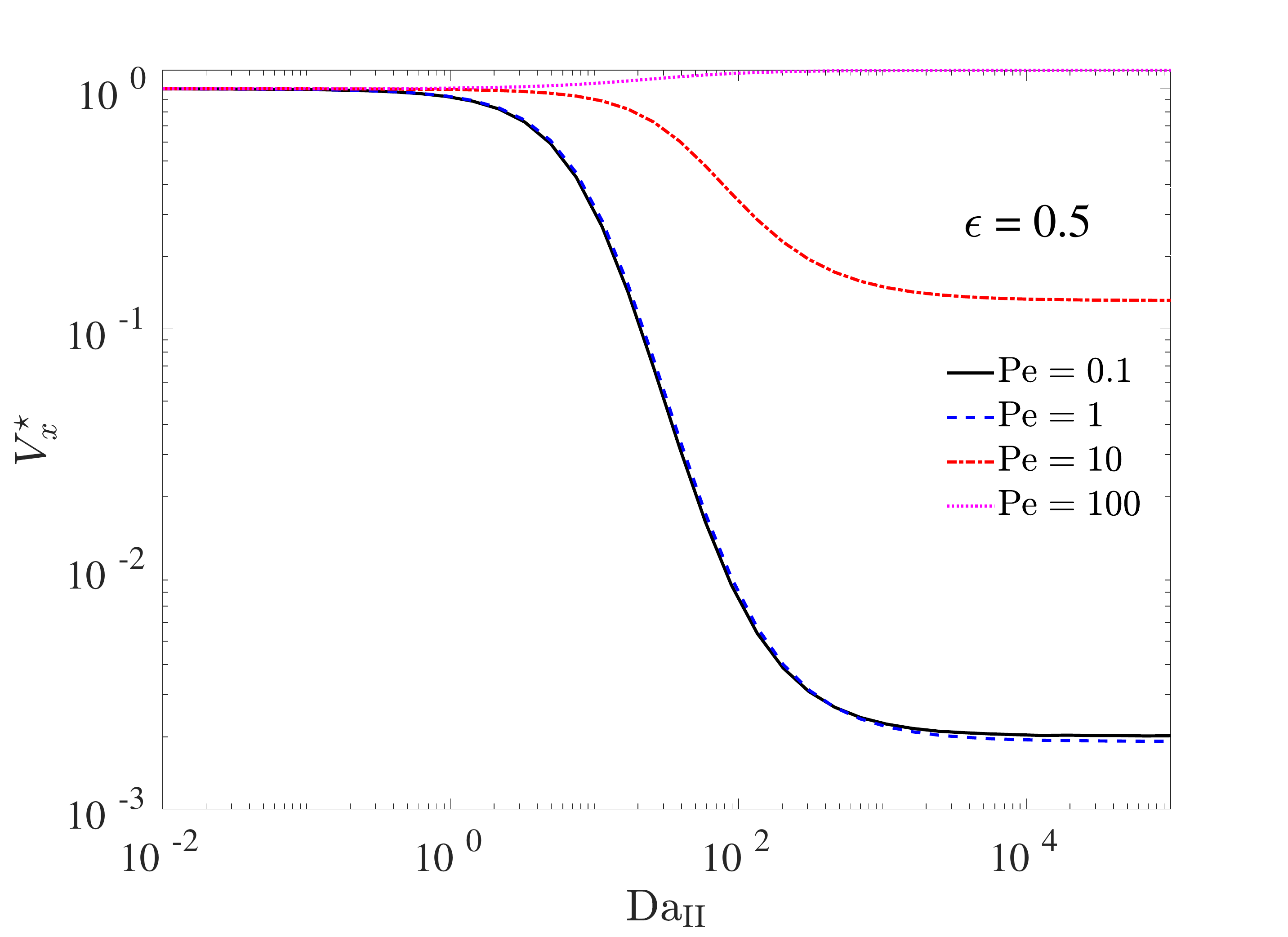}
        \caption{}
        \label{subfig:eps05_U}
    \end{subfigure}%
    \caption{Scaled effective parameters as a function of the microscopic parameters for $\epsilon=0.5$.}
    \label{fig::eps05}
\end{figure}

We illustrate how the present method can be applied to large scale studies by studying the effect of $\PE$ and $\DAII$ on the effective parameters of the homogenised transport equation.
Results are presented in Figures~\ref{fig::eps01},~\ref{fig::eps07} and~\ref{fig::eps05} for different values of the porosity $\eps$. Furthermore, full numerical results are provided in the additional material.

We probe the range  $ \DAII \in \left[ 10^{-2}, 10^{5} \right]$ in order to capture both the limits tending to Neumann and Dirichlet boundary conditions. 

Results show that $\lambda$ increases with $\DAII$ in an exponential fashion and reaches a saturation value around $\DAII=100$. This is consistent with other studies in literature \citep{boccardo2018robust,Bourbatache2020}. Furthermore and as expected,  $\lambda$ increases with decreasing  $\epsilon$ since the exchange surface per unit volume increases. Increasing $\PE$ also leads to an increased interphase transfer, as documented, for example, in a number of correlations for dense particle beds \cite{Gunn1978,Deen2014,Municchi2017}.

The effective dispersion also shows a consistent behaviour for all values of $\epsilon$, and its dependence on $\DAII$ is in agreement with results in literature for diffusion--reaction. The increase with smaller $\epsilon$ can be explained considering that the flow field becomes more tortuous when larger portions of the domain are occupied by the solid phase. Interestingly, we notice that the effective dispersion coefficient in reactive systems is generally smaller with respect to the non-reactive case, with one noticeable exception around $\PE=10$ (see panel c of Figures \ref{fig::eps01} and \ref{fig::eps07}) for the transverse diffusion. This may hint the existence of a region in the parameter space where the effective dispersion behaves irregularly. However, in the current work, we do not have enough samples to investigate such a feature.  

It is interesting to observe that the effective dispersion increases with increasing P\'eclet number for $\DAII>1$. The fact that the effective dispersion increases with increasing $\PE$ is well documented in literature  \cite{Delgado2006}, and correlations for random arrays of spheres are able to reproduce results from direct numerical simulations with reasonable accuracy \cite{Municchi2017}. Results suggest that advection is capable of mitigating the reduction of effective dispersion induced by reactions, and that it can even reverse the trend, leading to an enhancement of the effective dispersion.

Results from different values of $\epsilon$ show little variation in the profiles of $\lambda$, the values of which increases for higher P\'eclet and $\epsilon$, while showing an interesting behaviour for both the effective dispersion and velocity.
Specifically, the effective dispersion profiles do not show any consistent trend with the P\'eclet number, which indicates the complex interconnection between this and the Damk\"ohler number on the effective transport properties (through the corrector $\boldsymbol{\chi}$). However, it can be seen that both $\mathcal{D}'_{xx}$ and $\mathcal{D}'_{yy}$  tend to occupy a larger range of values with decreasing porosity. This can be explained considering that the flow field becomes more tortuous when larger portions of the domain are occupied by the solid phase. Interestingly, we notice that the effective diffusion coefficient in reactive systems can be both larger or smaller with respect to the non-reactive case. At low P\'eclet numbers the effective dispersion approaches the profiles obtained in previous works on the homogenisation of the diffusion equation \citep{Bourbatache2020}. It is interesting to observe that the effective dispersion increases with increasing P\'eclet number for $\DAII>1$. The fact that the effective dispersion increases with increasing $\PE$ is well documented in literature  \cite{Delgado2006}, and correlations for random arrays of spheres are able to reproduce results from direct numerical simulations with reasonable accuracy \cite{Municchi2017}. Results suggest that advection is capable of mitigating the reduction of effective dispersion induced by reactions, and that it can even reverse the trend, leading to an enhancement of the effective dispersion. To understand the physical reason for this, one should consider three points:

\begin{itemize}
    \item[i] The Damk\"oler number can be understood as a parameter allowing to switch from a Neumann boundary condition (no gradient) to a Dirichlet boundary condition (largest gradient) in a smooth manner. Thus, the Damk\"oler number generates gradients.
    \item[ii] Increasing the P\'eclet number means thinning the boundary layer, thus increasing local gradients. Most importantly, high P\'eclet numbers induce wakes behind the cylinders that stretch the boundary layer in the flow direction, as observed in Figure~\ref{subfig::phi}. This gives rise to stronger gradients and, in extreme cases, even filaments. Thus the P\'eclet number propagates the gradients in the downstream direction.
    \item[iii] The effective dispersion is a function of the gradient of $\boldsymbol{\chi}$ (see Equation~\ref{eq::Deff_2} ), which is a measure of the local variation of $c$. 
\end{itemize}

Therefore, $\boldsymbol{\chi}$ depends on the Damk\"oler number through $\beta$ and the gradients of $\phi$ and $\phi^{\dagger}$ (see Equation~\ref{eq::first_order_corrector_eq}). Notice that in the case of pure diffusion one would have a self-adjoint problem $\phi_{sa} = \phi^{\dagger}_{sa}$ (the subscript $sa$ indicates that they solve the self adjoint problem) and thus $\bv^{\star}=0$. However, when advection is present one has $\phi \neq \phi^{\dagger}$, and $\ny (\phi-\phi_{sa})$ is generally similar in magnitude but opposite in sign to $\ny (\phi^{\dagger}-\phi_{sa}^{\dagger})$. Due to the symmetry of the geometry we used, they are actually of the same magnitude but opposite sign (see Figure~\ref{subfig::phi} and Figure~\ref{subfig::phiAdj}). This implies that the strength of the advection term  in Equation~\ref{eq::first_order_corrector_eq}, which is given by $\bv^{\star}$, is not only affected by the P\'eclet number, but also by the Damk\"oler number. However, when $\PE$ increases, this effect becomes less important and the effective dispersion becomes almost independent on $\DAII$.

The effective velocity is probably the macroscopic parameters which exhibits the most interesting behaviour. In fact, for low $\PE$ it increases and then decreases with $\DAII$, while it is an increasing function of $\DAII$ for high $\PE$. This effect is more pronounced for high $\eps$. This phenomenon reveals an interesting role of $\beta$, which acts as a weighting function of the velocity field. Since $\beta$ is larger inside the channel than at the boundaries, the regions at higher velocity are weighted more than the boundary layers. This results in a channelling effect, leading to an apparent velocity higher than the average velocity. Notice that such effect is present in a reduced manner at lower $\epsilon$. This happens because $c$ is consumed by the surface reaction, resulting in a slowing down of the advective flux.  

Although these two dimensional results are not directly applicable to many real porous media, these results lead to a better understanding of reactive flow through ordered arrays of cylinders and may provide useful hints for design of membranes, heat exchangers or catalysts.

\subsection{Deposition of charged colloids}
\label{SS:colloids}
The method presented in this work allows us to find reliable macroscopic models for a wide range of other problems that present similar mathematical difficulties. As an example, we present here results for the upscaling of  colloidal transport in presence of a Coulomb-like attraction potential. Similarly,  other more realistic (e.g., DLVO \citep{Ohshima2014} and gravity) potentials can be easily implemented in the mathematical and numerical method.

We consider here a multiple inverse--distance potential (here presented in dimensionless form):
\begin{equation}
    \label{eq::pot}
    \Lambda \ofy = \sum \limits_{i=0}^{N_b}\frac{ \Lambda_i}{r_i\ofy + R_i} \, ,
\end{equation}
where the summation is carried over the bodies $i$ that generate the field and $N_b$ is the number of bodies generating a potential field. We denote the magnitude of the distance from the surface of body $i$ as $r_i$ and $R_i$ is a reference (dimensionless) length (e.g., the radius in case of spherical bodies). Each body is characterised by a potential strength $\Lambda_i$. As we have already seen, 
potential forces in the Smoluchowski approximation result in  a net velocity $\bm{v}_{\Lambda} \ofy$, which is generally not solenoidal:
\begin{equation}
    \label{eq::potVel}
  \bm{v}_{\Lambda} \ofy =  - \frac{1}{\eps}\ny \Lambda \ofy = -\sum \limits_{i=0}^{N_b}\frac{\Lambda_0}{\eps} \ny \left(\frac{1}{r_i\ofy + R_i}\right) \, ,
\end{equation}
where, for simplicity, we set a constant $\Lambda_i=\Lambda_0$, so that each grain has the same potential strength, and choose $R_i$ simply as the dimensionless grain radii.
We present here results for $\PE=1$ and with a ratio between potential and advective terms $\mu=\frac{\Lambda_0}{\PE}$ ranging from $-10$ (attractive potential) to $10$ (repulsive potential). We limit the study here to the case of $\DA=0$.
We remind here that the total velocity for this case is not solenoidal and non-zero at the wall. This means that, even without reactions, the standard homogenisation would fail as the boundary conditions are of mixed (Robin) type.

\begin{figure}[htbp]

    \centering
    \begin{subfigure}[htbp]{0.5\textwidth}
        \centering
        \includegraphics[width=0.87\textwidth]{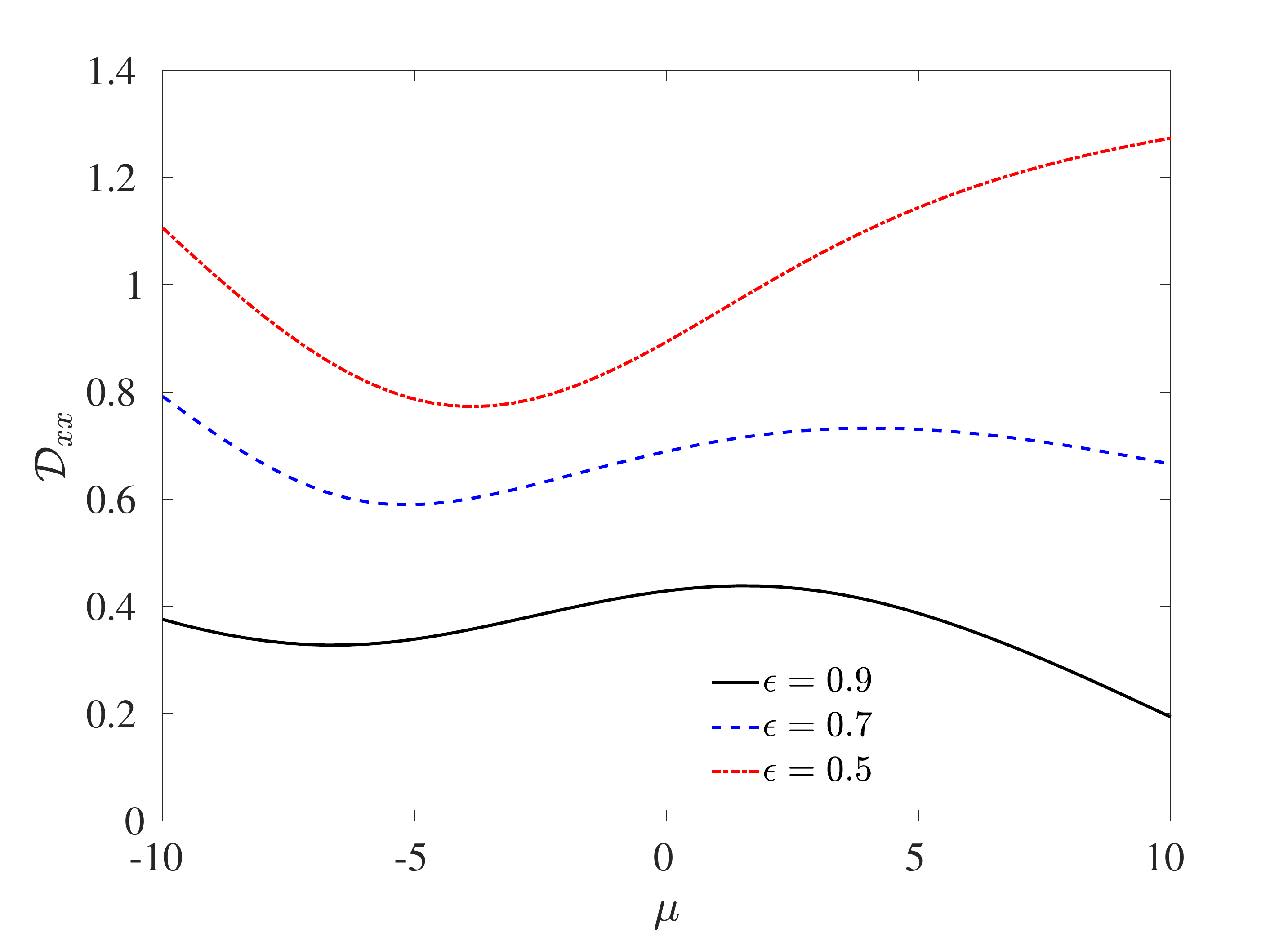}
        \caption{}
        \label{subfig::rep_Dx}
    \end{subfigure}%
    ~ 
    \begin{subfigure}[htbp]{0.5\textwidth}
        \centering
        \includegraphics[width=0.9\textwidth]{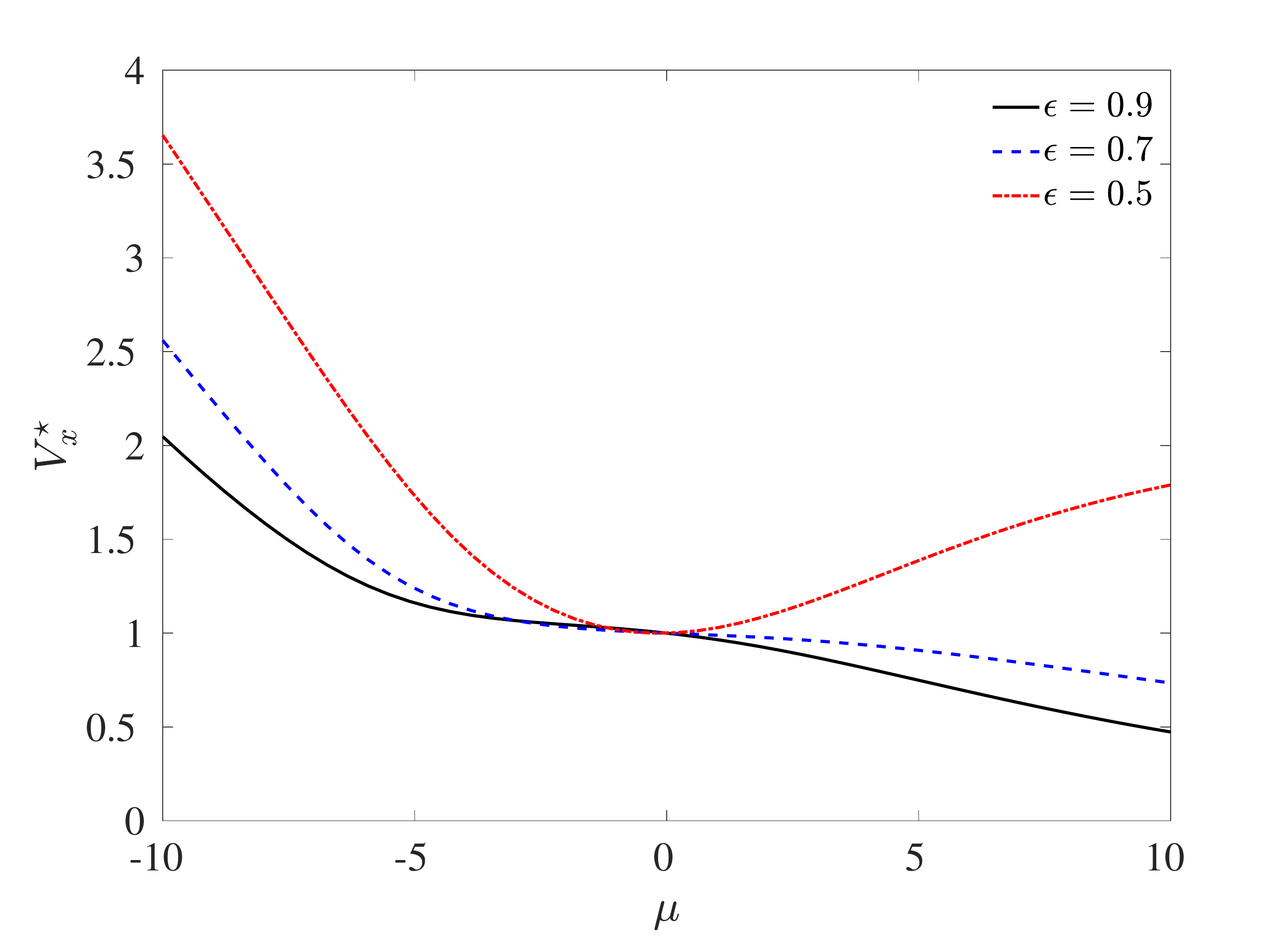}
        \caption{}
        \label{subfig::rep_Vx}
    \end{subfigure}%
    \caption{Scaled effective parameters as a function of  the potential number $\mu$ for different values of the porosity $\eps$. $\PE$ is fixed to $1$ while $\DAII$ is zero. Notice that a negative value $\mu$ correspond to a repulsive potential, while a positive value corresponds to an attractive potential.}
    \label{fig::coll}
\end{figure}

In Fig.~\ref{fig::coll} we show results for $\PE=1$ and $\DAII=0$ (non reactive) over a range of values of $\mu$ and $\epsilon$. When the potential is attractive (positive $\mu$) the effective dispersion increases with $\mu$ while a repulsive potential leads to an overall decrease in $\mathcal{D}_xx$ followed by a further increase with a minimum around $\mu=-6$. This is expected, since the repulsive potential will tend to scatter the particles far from the inclusions. Similarly, attractive potentials will drive the particles into regions of low local velocity (i.e., close to the boundary layers), while repulsive potentials will push the particles towards regions of flow channelling. This phenomenon is expected to have a stronger effect for low values of the porosity, where the channelling is stronger. However, since the magnitude of the velocity field is lower at the centre of the cell than in the regions immediately above and below (see Fig.\ref{subfig::uplus}), even an attractive potential can  result in an increase of the effective velocity, especially for low values of the porosity, where the channelling is stronger. This fact can also explain the non-uniform dependence of $\mathcal{D}_{xx}$ on $\mu$.

\subsection{Application to three-dimensional geometries}
\begin{figure}[htbp]

    \centering
    \begin{subfigure}[htbp]{0.5\textwidth}
        \centering
        \includegraphics[width=0.95\textwidth]{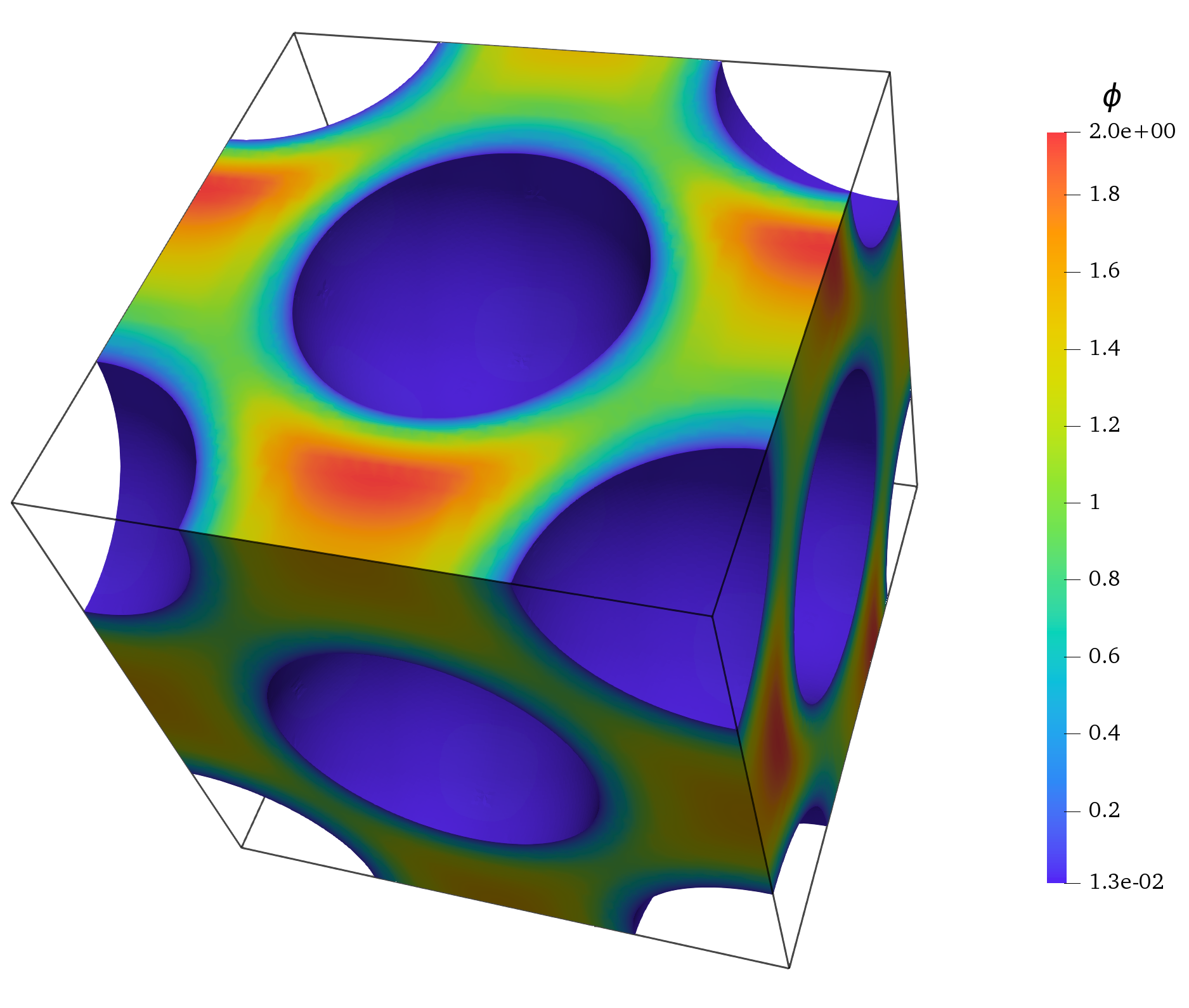}
        \caption{}
        \label{subfig::phi3d}
    \end{subfigure}%
    ~ 
    \begin{subfigure}[htbp]{0.5\textwidth}
        \centering
        \includegraphics[width=0.95\textwidth]{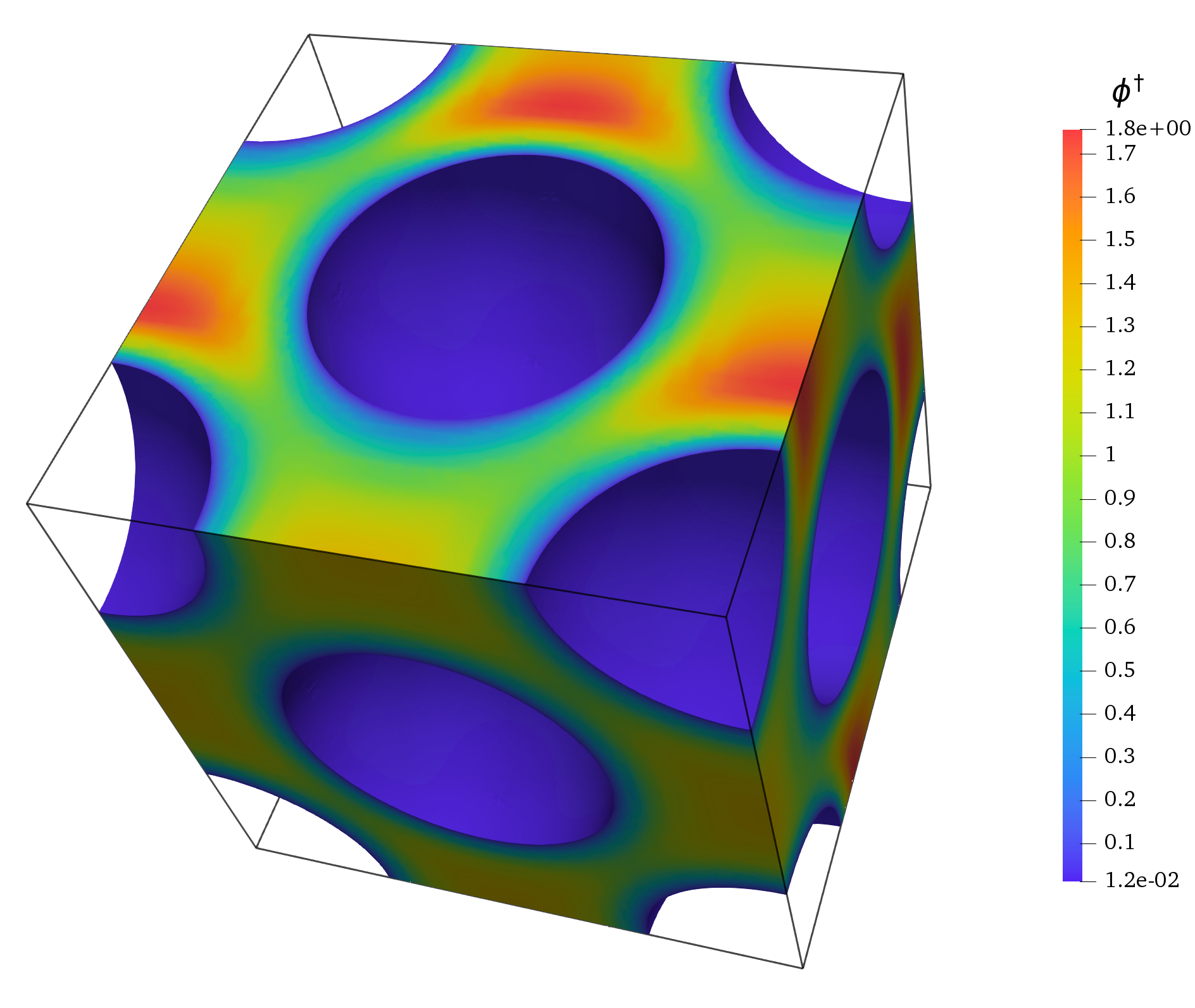}
        \caption{}
        \label{subfig::phiadj3d}
    \end{subfigure}%
    \caption{Contour plots of $\phi$ and $\phi^{\dagger}$.}
    \label{fig::3dpsi}
\end{figure}

\begin{figure}[htbp]

    \centering
    \begin{subfigure}[htbp]{0.5\textwidth}
        \centering
        \includegraphics[width=0.95\textwidth]{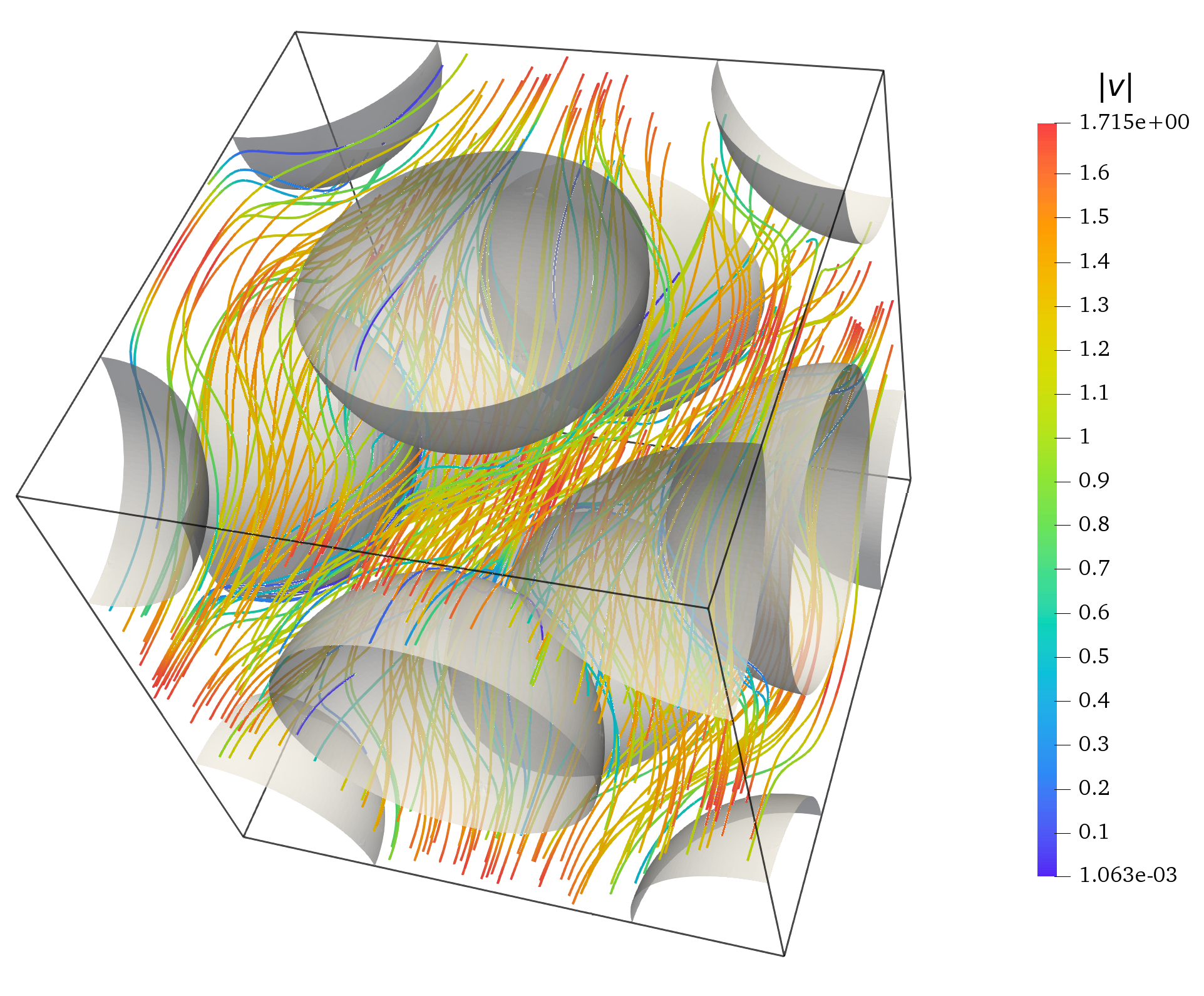}
        \caption{}
        \label{subfig::u3d}
    \end{subfigure}%
    ~ 
    \begin{subfigure}[htbp]{0.5\textwidth}
        \centering
        \includegraphics[width=0.95\textwidth]{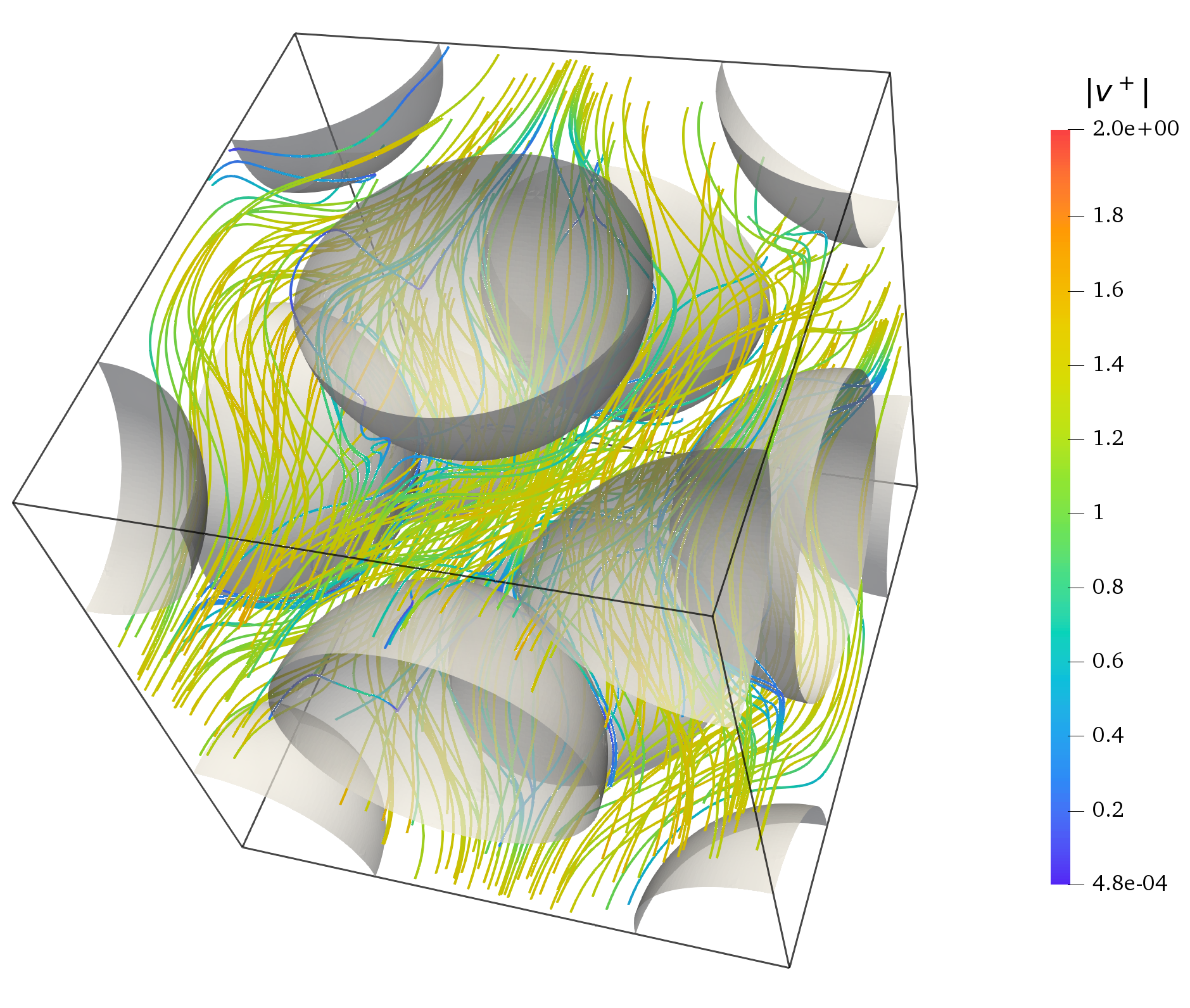}
        \caption{}
        \label{subfig::ulpus3d}
    \end{subfigure}%
    \caption{Comparison between streamlines of $\bv$ and $\bv^{+}$.}
    \label{fig::3du}
\end{figure}

\begin{figure}[htbp]
    \centering
    \includegraphics[width=0.5\textwidth]{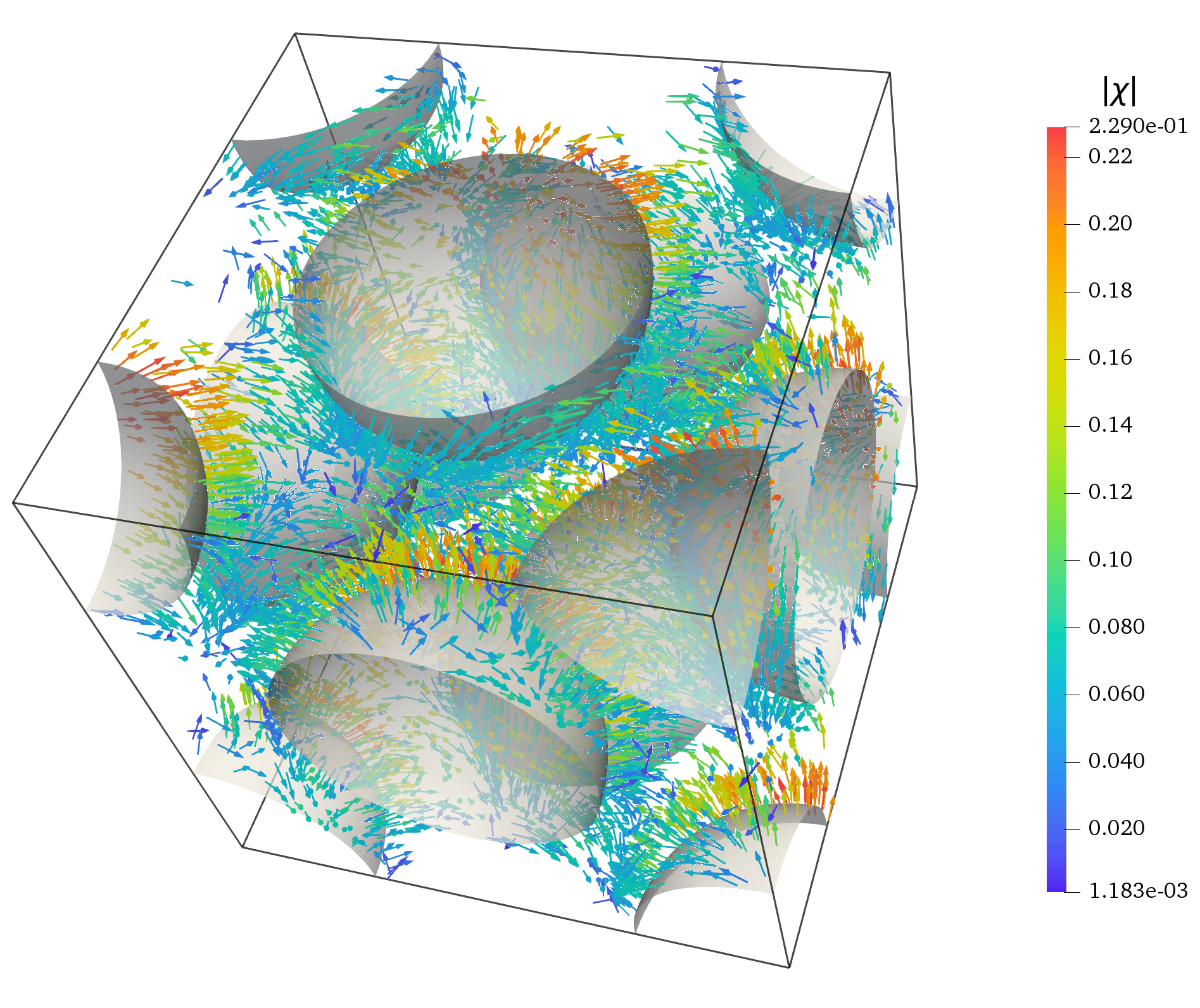}
    \caption{Vector field representation of $\boldsymbol{\chi}$.}
    \label{fig::3dchi}
\end{figure}

The numerical method presented in this work can be easily extended to three-dimensional geometries thanks to the flexibility of the OpenFOAM\textsuperscript{\textregistered} library. While a parametric study is beyond the scopes of the present work, we illustrate some results from three dimensional geometries in Figures~\ref{fig::3dpsi}, \ref{fig::3du}, and \ref{fig::3dchi}. In these simulations we employed $\DAII=1060$, $\PE=10$, and $\epsilon=0.9$.

It is interesting to observe how the effective parameters in the three dimensional configuration differ from those obtained from a two dimensional domain with the same microscopic parameters (i.e., with the same $\DAII$, $\PE$, and $\epsilon$). Table \ref{tab::3d2d} reveals that the eigenvalue is definitely the parameter that changes the most. This can be attributed to the larger value of the specific surface in the three-dimensional configuration, resulting in a larger exchange area in the unit cell and thus, in a larger $\lambda$. Another parameter increasing significantly is the longitudinal effective dispersion $\mathcal{D}^{\text{eff}}_{xx}$, which is a consequence of the increased tortuosity of the flow. No measurable change was observed for the transverse effective dispersion $\mathcal{D}^{\text{eff}}_{yy}$, which maintains the same value in the three-dimensional configuration. Finally, since the fluid is less constricted (being able to flow in one additional dimension), the effective velocity is lower in three dimensions, but still larger than the average velocity, thus hinting to a similar trend as in two dimensions.  

\begin{table}
\begin{center}
    \begin{tabular}{c|c|c}
         & 3D & 2D  \\
        \hline
        $\lambda$ & 5.346 & 1.523 \\
        $\mathcal{D}^{\text{eff}}_{xx}$ & 0.045 & 0.035 \\
        $\mathcal{D}^{\text{eff}}_{yy}$ & 0.071 & 0.071 \\
        $V^{\star}_{x}$ & 1.115 &  1.123\\
    \end{tabular}
    \caption{Comparison between effective parameters from a three-dimensional geometry and a two-dimensional geometry using the same microscopic parameters.}
    \label{tab::3d2d}
\end{center}
\end{table}

\section{Conclusions}
\label{S:conc}

In this work, we presented a methodology for the upscaling of reactive transport in porous media based on the works of \citet{ALLAIRE2007523} and \citet{Mauri_1991}. Such upscaling procedure has been described in details for homogeneous boundary conditions and extended to non-homogeneous boundary conditions and potential forces. This allows us to significantly extend the range of physical problems for which upscaled equations can be found, bypassing the usual limitation of slow reaction and divergence-free velocity fields.

Furthermore, to make the approach available to the community, we implemented the method in the open-source library OpenFOAM\textsuperscript{\textregistered} \cite{Foundation2014} and compared its predictions against fully resolved microscale simulations finding excellent agreement. This confirms the power and accuracy of this homogenisation-based approach to extend the applicability of macroscopic transport theory beyond the  'standard' problems that rely on solenoidal velocity field and slow reactions. It is important to notice that our approach does not include 'conjugate' transfer (i.e., transport inside the solid grain) but the present approach can be conveniently coupled with the recently proposed generalised multi-rate transfer model \cite{Municchi2020}.

Finally, we presented a parametric study of reactive transport in ordered arrays of cylinders to illustrate the usage of the proposed method and numerical code. We found that all the effective parameters (effective velocity, dispersion, and reaction) depend on the P\'eclet and Damk\"ohler number in a complex manner and that the method is able to correctly recover the limiting cases of Dirichlet (infinitely fast reaction) and Neumann (non-reactive) boundary conditions. This study is extended for a case of charged colloids driven by a flow field and an attractive electrostatic force between the particles and the solid grains. The additional drift term in the equation is not divergence free. This means that standard homogenisation cannot be applied. For this test case, all the effective parameters shows a non-trivial non-linear dependence on the potential.

While we restricted our numerical investigation to FCC configurations, the presented methodology and numerical code can be seamlessly extended to any other geometry. In future works, we plan to use this method to study heterogeneous bulk reactions also at the macroscale, as well as heat/mass transfer in grain packings and suspensions \cite{Municchi2017,Municchi2018}. Since the method relies only on stationary cell problems, it is also particularly suitable for computing statistics of random porous media \cite{icardi2016predictivity}. In fact, when computing the effective parameters of a random system one is often interested in the asymptotic regime, at which such parameters are constant and bear no dependence on the initial conditions. Conventionally, the strategies employed to achieve the asymptotic regime consist in employing some pseudo-periodic boundary conditions \citep{Municchi2017,Tenneti2013} or in constructing large domains, assuming that the asymptotic regime is established far from the boundaries \citep{Tavassoli2016,Deen2014}. Another approach consists in employing standard periodic boundary condition and letting the system evolve towards a state of saturation (constant concentration everywhere), until the effective parameters are not depending on time anymore \citep{Derksen2014}. On the contrary, the method we propose is based on the solution of cell problems that do not depend on time and are valid in the asymptotic regime. Thus, it requires simple periodic boundary conditions and no time-dependent equations need to be solved, removing any dependence on the initial conditions. This makes the method very efficient for large studies involving a wide range of parameters.

The study of charged colloids can be extended considering other forces (gravity) and more complex, non-linear and electrokinetic boundary conditions \cite{joekar2019coupled} with important applications for electrokinetic energy conversion in nanofluidic channels \cite{Ren2008}.

As a final remark, our numerical code is open-source and freely available \cite{Municchi2020b} with the objective of extending the use of homogenisation-based techniques to a wider community and providing an 'upscaling toolbox' with solid mathematical foundations.

\section*{Acknowledgements}
This work has been funded by the European Union's Horizon 2020
research and innovation programme, grant agreement number 764531, "SECURe -- Subsurface Evaluation of Carbon capture and storage and Unconventional risks".


\clearpage
\bibliographystyle{elsarticle-num-names} 
\bibliography{references}


\end{document}